\documentclass{jnmp}

\usepackage{amsmath}
\usepackage{epsfig}
\JNMPnumberwithin{equation}{section}

\theoremstyle{definition}
\newtheorem*{remark}{Remark}

\newcommand{\ddt}{\partial \over \partial t}
\newcommand{\ddx}{\partial \over \partial x}
\newcommand{\ddu}{\partial \over \partial u}
\newcommand{\dds}{\partial \over \partial s}

\def\dsd#1{ \mathop{#1}\limits_{+s}}

\def\da#1{ \mathop{#1}\limits_{+\tau}}
\def\db#1{ \mathop{#1}\limits_{-\tau}}
\def\dc#1{ \mathop{#1}\limits_{\pm \tau}}
\def\dd#1{ \mathop{#1}\limits_{+h}}
\def\df#1{ \mathop{#1}\limits_{-h}}

\def\dg#1{ \mathop{#1}\limits_{\pm h}}

\def\dh#1{ \mathop{#1}\limits_ h}

\begin{document}
\renewcommand{\evenhead}{V Dorodnitsyn and R Kozlov}
\renewcommand{\oddhead}{A Heat Transfer with a Source:
the Complete Set of Invariant Difference Schemes}

\thispagestyle{empty}

\FirstPageHead{10}{1}{2003}{\pageref{dorodnitsyn-firstpage}--\pageref{dorodnitsyn-lastpage}}{Article}

\copyrightnote{2003}{V Dorodnitsyn and R Kozlov}

\Name{A Heat Transfer with a Source: \\
the Complete Set of Invariant Difference Schemes}
\label{dorodnitsyn-firstpage}

\Author{Vladimir DORODNITSYN~$^\dag$ and Roman KOZLOV~$^\ddag$}

\Address{$^\dag$~Keldysh Institute of Applied Mathematics,
Miusskaya Pl.~4,     Moscow, 125047, Russia \\
~~E-mail: dorod@spp.Keldysh.ru\\[10pt]
$^\ddag$~Department of Informatics, University of Oslo, 0371, Oslo, Norway \\
~~E-mail: kozlov@ifi.uio.no}

\Date{Received February 28, 2002; Revised July 11, 2002;
Accepted July 24, 2002}

\begin{abstract}
\noindent
In this letter we present the set of invariant
difference equations
and meshes which preserve the Lie group symmetries of
the equation
$u_{t}  = ( K(u) u_{x} )_{x} + Q(u)$.
All special cases of $K(u)$ and $Q(u)$
that extend the symmetry group admitted by the
differential equation are considered.
This paper completes the paper [{\it J. Phys. A: Math. Gen. } { \bf 30},
Nr.~23 (1997),  8139--8155], where a few invariant
models for heat transfer equations were presented.
\end{abstract}

\section{Introduction}

Symmetries are fundamental features of the differential
equations of mathematical physics. It yield a number of
useful properties such as  integrability of
ODEs, symmetry reduction of PDEs,
existence of various types invariant solutions,
conservation laws for the invariant variational problems etc.
Therefore, preserving symmetries in discrete schemes, we retain
qualitative properties of the underlying differential equations.

The purpose of this paper is to develop the entire set of invariant
difference schemes for the heat transfer equation
\begin{equation} \label{eq0}
u_t = \left( K(u) u_{x} \right) _{x} + Q(u),
\end{equation}
for all special cases of the coefficients $K(u)$, $Q(u)$ which extend the
symmetry group admitted by equation (\ref{eq0}).
This paper is based on the Lie group classification~\cite{[4]}
(see also~[1]) of the equation (\ref{eq0})
with arbitrary $K(u)$ and $Q(u)$. This classification
contains the result of L~V~Ovsyannikov \cite{[18]} for equation~(\ref{eq0})
with $Q \equiv 0$ as well as symmetries of the linear case $K \equiv 1$,
$Q \equiv 0$, which were known by S~Lie.

A few examples of the invariant difference schemes and meshes were
considered in~\cite{[Bakir]}. In the present
paper we complete the paper~[2] going through all cases of $K(u)$ and $Q(u)$
identified in the group classification~\cite{[4]}, we construct
difference equations and meshes (lattices) which admit the same Lie groups
of point transformations as their continuous limits.

Lie group analysis of difference equations is a very active field
of research where many contributions were done, and various
approaches were applied by several authors (see~[18]).
In our approach which we are following in this paper we pose the question:
How does one discretize a differential equation while preserving all of its Lie
point symmetries? Thus a~differential equation and its Lie group symmetry
are {\it a priory} given but not a difference model. One then looks for a difference scheme,
i.e.\ a difference equation and a mesh, that have the same symmetry group
and the same Lie algebra. The basic steps in this direction were done [5, 6, 7, 8, 9, 10, 11],
 which were summarized in a recent book~\cite{[12]}.
The main idea is that the invariant difference equations
and meshes can be constructed with the help of the entire set of difference invariants of the
corresponding Lie group. In the next section we explain how to construct
difference models that conserve the whole group of point transformations admitted
by the differential equations.

The article is organized as follows. Section~2
provides a brief overview of the invariant discretization procedure. In
Sections~3, 4, 5 and~6 we consider the cases of
an arbitrary heat transfer coefficient $K(u)$,
the exponential heat transfer coefficient $e^u$,
the power heat transfer coefficient $u^{ \sigma }$ and
the special case of power heat transfer coefficient:
$u^{-4/3}$
corres\-pondingly.

  Section~7 is devoted to the linear heat conductivity with
a source. In particular, this section covers detailed
study of the invariant difference scheme for the linear heat equation without
a source ($Q = 0$) including such aspects as superposition principle, reduction
of the invariant scheme on the optimal system of subalgebras and the
way to transform the moving mesh scheme into a stationary one.
Notice that in the paper~[19] there were considered some difference
approximations of the linear heat transfer equation, which preserve its different
symmetries on different meshes. In [19] a difference equation and
 a~mesh are {\it a priori} given, then it was shown that for some
kind of the mesh there were preserved some symmetries of the linear
heat equation and another meshes preserve other parts of the
symmetries. 
Thus, there are no difference schemes which conserve the {\it entire set of symmetries}
in one difference model. In Section~7 we will develop the difference mesh
and difference equation, which conserve the complete set of original symmetries
in the one and the same difference scheme.

The same approach we will apply for some tens of other nonlinear models
of heat equation (\ref{eq0}).
Summarizing conclusions end up the consideration of  the entire
set of invariant schemes for the equation (\ref{eq0}).

\section{Symmetry preserving discretization procedure}

1. Let us briefly describe this method, which was called
{\it the method of finite-difference invariants}~\cite{[6]}.

Let the differential equation
\begin{equation} \label{eqet}
F(t,x,u ,u_t, u_x, \ldots) = 0
\end{equation}
admit a known symmetry group $G_n$, whose Lie algebra
is spanned by the
operators $X_1, \ldots,X_n$ of the form
\begin{equation} \label{operator1}
X _i = \xi  ^{t} _i {\partial \over \partial t} +
\xi ^{x} _i {\partial \over \partial x} +
\eta _i   {\partial \over \partial u} , \qquad i = 1,
\ldots, n,
\end{equation}
where the coefficients
$ \xi ^{x} $, $ \xi ^{t}  $ and $\eta $ are functions
of $t$, $x$, $u$, since we consider Lie point symmetries.

Then we would like to propose a discrete model
\begin{gather}
F(z) = 0 ,\nonumber\\
\Omega (z,h) = 0 ,
\label{eqtt}
\end{gather}
where the first equation is the approximation of the initial
differential equation  and the second one defines a difference mesh.
Both of these equations {\it a priori} are not given and we
have to establish  the invariant mesh, on which we should
approximate the original heat equation.
We will show, that in all special cases of equation~(\ref{eq0})
one can construct the system (\ref{eqtt}),
starting from the entire set of finite-difference invariants of the corresponding Lie group.

 We denote by $z$ in (\ref{eqtt}) a finite number of difference variables,
which are used in the considered {\it difference stencil}, i.e.\ a finite number
of mesh points, which are needed for the approximation of the differential
equation (\ref{eqet}). The equations (\ref{eqtt})
can be explicitly connected with each other (if invariant mesh depends on solution) or not.
In the last case we can choose the invariant mesh firstly (for example, a fixed mesh)
and then construct the invariant approximation of the original
equation. 
If a mesh depends on the solution, all specifications of the mesh made
in advance lead to
restrictions on the symmetries which may be admitted by the considered
discrete models.

2. The idea of the method of finite-difference invariants springs from
the invariant representation of differential equations. In the continuous case for  the group
$G_n$ we can find the complete set of
functionally independent differential invariants $ J = ( J_1,  J_2, \ldots, J_k) $
in the
specified space which contains dependent and
independent variables as well as
the set of derivatives up to the highest derivatives
involved in
the formulation of the PDE~\cite{[16]}.
For the heat equation~(\ref{eq0}) we
consider the space $M \sim ( t,x,u,u_t,u_x, u_{xx})$.
We prolong the operator (\ref{operator1}) on the
variables of the
space~$M$
\[
{\bf pr} \; X  = X  +
{\zeta}^t  {\partial \over \partial u_t} +
{\zeta}^x {\partial \over \partial u_x} +
{\zeta}^{xx}  {\partial \over \partial u_{xx} },
\]
with
\begin{gather*}
{\zeta}^t  = D_t (\eta) - u_t D_t (\xi ^t ) - u_x D_t
(\xi ^x ),
\qquad
{\zeta}^x  = D_x (\eta) - u_t D_x (\xi ^t ) - u_x D_x
(\xi ^x ),\\
{\zeta}^{xx}  = D_x ({\zeta}^x ) - u_{tx} D_x (\xi ^t)
- u_{xx} D_x (\xi ^x ),
\end{gather*}
where $D_t$ and $D_x$ are the total derivative
operators for time
and space correspondingly.
Differential invariants are solutions of the system of linear equations
\[
  {\bf pr} \; X_i \; \Phi (t,x,u,u_t,u_x, u_{xx}) = 0 ,
\qquad i = 1, \ldots, n,
\]
and can be solved by standard procedure (see~[15]).

Then we represent the invariant differential equation
in terms of
these invariants
\[
\tilde{F} ( J_1,  J_2, \ldots, J_k) = 0 .
\]
The obtained equation is invariant with
respect to
the group $G_n$.

3. In the discrete case the situation is more complicated. Any given differential
equation can be approximated by means of infinitely many difference equations
and meshes, which have the original differential equation as its continuous limit.
The requirement of preservation of Lie group properties of the differential
equation in its discrete counterpart still leaves some freedom in the approximations.
Thus, as it can be seen recently, at some point we have to make a chose among
the general family of invariant meshes.

 The structure of the admitted group essentially effects on construction
of equations and meshes. Group transformations can break the
geometric structure of the difference mesh that influences approximation and
other properties of difference equations. First steps to the
construction of the difference grids geometry based on the symmetries of
the initial difference model were done in~[5, 6, 7, 9]. There were found classes
of transformations that conserve uniformity, orthogonality and other
properties of the grids.

It was shown [5, 6, 7, 9] that a transformation defined by
(\ref{operator1})
conserves uniformity of a mesh in $t$ and $x$ directions if and only if
\begin{gather} \label{ct}
{ \da{D} }{ \db{D}  }( \xi ^{t} ) = 0 ,\\
\label{ch}
{ \dd{D} }{ \df{D} }  ( \xi ^{x} ) = 0 ,
\end{gather}
where ${\dc{D}}$ and ${\dg{D}}$ denote total
difference
derivatives in the time and space directions
with steps~$\tau$ and~$h$ correspondingly.

For an orthogonal mesh to be conserved under the
transformation, it is necessary and sufficient that
\begin{equation} \label{cht}
{\dd{D} }( \xi ^{t} ) = -{ \da{D}} ( \xi ^{x} ) .
\end{equation}

When condition (\ref{cht}) is not satisfied for a
given group, the
flatness of the layer of a~grid in some direction is
rather important. For evolution equations it is significant to have flat
time layers.
There is a simple
criterion  of the invariance of flat time layers
under the action of a given operator~(\ref{operator1}):
\begin{equation} \label{cht2}
{\dg{D}} { \da{D}} ( \xi ^{t} ) = 0.
\end{equation}
These condition
specify
invariant geometry of grids
for the given Lie group symmetries.

If the operator coefficients $\xi ^{t}$, $\xi ^{x}$ do not depend
on solution, then we can choose the invariant mesh as any solution
of corresponding condition (\ref{ct})--(\ref{cht2}). Otherwise the
conditions (\ref{ct})--(\ref{cht2}) should hold on the solutions
of the considered difference model. In that case we can figure a
mesh out starting from the set of difference invariants.

Further we choose a stencil
which is sufficient to approximate all derivatives
which appear in the equation. We will consider
six-point stencils
which have three points on each of two time layers.
Such stencils
allow us to write down both explicit and implicit
difference
schemes. For different transformation groups we will
consider
different meshes: orthogonal mesh which is uniform in
space,
orthogonal mesh which is nonuniform in space and
nonorthogonal in time-space mesh, i.e.\ moving mesh.
The corresponding stencils are
different. Furthermore,
the corresponding spaces of discrete variables are
of different dimensions so that they have
different number of difference invariants $ I = ( I_1,
 I_2, \ldots, I_l) $
for the same Lie group~$G_n$.

For example, let us take an orthogonal mesh which is
uniform in space.
(We will describe later for which groups this mesh can
be considered.)
The stencil of this mesh is shown in Fig.~1.
\begin{figure}[th]  \label{grort}
{
\vspace*{-5mm}
\begin{picture}(300,100)
\put(105,0){\begin{picture}(200,100)%
\put(30,30){\line(1,0){140}}
\put(30,70){\line(1,0){140}}
\put(100,30){\line(0,1){40}}
\put(30,30){\circle*{5}}
\put(30,70){\circle*{5}}
\put(170,70){\circle*{5}}
\put(170,30){\circle*{5}}
\put(100,30){\circle*{5}}
\put(100,70){\circle*{5}}
\put(5,75){$(x-h,\hat{t},\hat{u}_{-})$}
\put(85,75){$(x,\hat{t},\hat{u})$}
\put(145,75){$(x+h,\hat{t},\hat{u}_{+})$}
\put(5,20){$(x-h,t,u_{-})$}
\put(85,20){$(x,t,u)$}
\put(145,20){$(x+h,t,u_{+})$}
\end{picture}}
\end{picture}}

\vspace{-5mm}

\caption{The stencil of the orthogonal mesh.}
\end{figure}
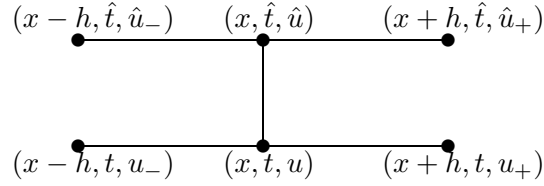

The corresponding discrete subspace is ten-dimensional:
$ M \sim
(t,x,\tau,h,u,u_{-},u_{+},\hat{u}$, $\hat{u}_{-},\hat{u}_{+})$,
where $  \tau = \hat{t} - t $. The prolonged operator
(\ref{operator1}) in this subspace has the form
\begin{gather*}
{\bf pr} \; X    = \xi  ^{t} {\partial \over \partial
t} +
  \xi  ^{x} {\partial \over \partial x} +
(   \hat{\xi} ^{t} -   \xi  ^{t}) {\partial \over
\partial \tau  } +
 ( \xi _+ ^{x} -  \xi ^{x} )  {\partial \over \partial
h}\\
\phantom{{\bf pr} \; X    =}{}
+ \eta   {\partial \over \partial u} +
\eta _- {\partial \over \partial u_-} +
\eta _+  {\partial \over \partial u_+} +
 \hat{\eta}   {\partial \over \partial  \hat{u} } +
 \hat{\eta} _-  {\partial \over \partial  \hat{u}_-} +
 \hat{\eta} _+  {\partial \over \partial  \hat{u}_+} ,
\end{gather*}
where we use time and space shifts notations
$\hat{f} =  f(t+ \tau,x,u)$, ${f}_{-} =
f(t,x-h,u)$,
${f}_{+} =  f(t,x+h,u)$. The number of
functionally independent
invariants is given by
\begin{equation} \label{cnd}
l = \dim M - \mbox{rank}\; Z, \qquad l \geq 0 ,
\end{equation}
 with $\dim M =10$ and
the matrix $Z$ composed by the coefficients of the
prolonged
on the space~$M$ operators
\[
Z = \left(
\begin{array}{cccccccc}
\xi_1^t &    \xi_1  ^{x} &  ( \hat{\xi} _1^x  - {\xi}
_1^x )   &
( ({\xi} _1^x)_+  - {\xi} _1^x ) &  {\eta} _1& \cdots &
(\hat{\eta} _1) _+ \\
\vdots \\
\xi_n^t &    \xi_n  ^{x} & (  \hat{\xi} _n^x  - {\xi}
_n^x  )  &
( ({\xi} _n^x)_+  - {\xi} _n^x ) &  {\eta} _n  & \cdots &
(\hat{\eta} _n) _+ \\
\end{array}
\right) .
\]

Having found the finite-difference invariants as the solutions of system of linear equations
\[
  {\bf pr} \; X_i \;   \Phi
(t,x,\tau,h,u,u_{-},u_{+},\hat{u},\hat{u}_{-},\hat{u}_{+})
  = 0 ,
\qquad i = 1, \ldots, n,
\]
we can use them to approximate the differential
invariants
\[
J_j  = f_j (  I_1,  I_2, \ldots, I_l ) + O ( \tau
^{\alpha} , h ^{\beta} ) ,
\qquad j = 1,\ldots,k,
\]
where $\alpha$ and $\beta$ define some fixed
order of approximation.
Notice, that approximation error $O ( \tau ^{\alpha} ,
h ^{\beta} )$ is
invariant together with other terms in the above
representation.
{\it Substitution of difference invariants $I_i$ instead
of differential ones $J_i$ into the function  $   \tilde{F} $
provides us with an invariant difference scheme.}
Practically  we can often omit the
representation of the differential equation
in terms of its invariants and just approximate
the original differential equation by the finite-difference
invariants. The use of finite difference invariants
is the main point in both ways.

So, the first step in the invariant approximation
is the choice of the invariant mesh. The last step is the choice of  the invariant
discretization of the original equation on the invariant mesh.

The described above method is algorithmic. We would like to stress
that the invariant approximation in our way is still not unique. For example
extending the stencil (means enlarging the number of mesh points involved
in approximation) we can find invariant approximations of any higher order.

\section{An arbitrary heat transfer coefficient $\boldsymbol{K(u)}$}

Now we start to develop invariant schemes going throw all cases of the Lie group
classification~[4]. Let us note that the group classification of the
equation~(\ref{eq0}) was
done in~\cite{[4]} (see also~\cite{[1]})  up to
equivalent transformations:
\begin{equation} \label{tr}
\bar{t} = at + e, \qquad \bar{x} = bx+ f, \qquad
\bar{u} =cu + g,
\qquad \bar{K} = { b^{2} \over a}\, K, \qquad \bar{Q} =
{ c \over  a} \,Q,
\end{equation}
where $a$, $b$, $c$, $e$, $f$ and $g$ are arbitrary
constants,
$abc  \neq 0$.
These transformations do not change differential
structure of the equation~(\ref{eq0}), transforming an admitted group into
a~similar group of point transformations.

{\bf 1.} We start from general case, when the coefficients
 $K(u)$ and $Q(u)$ are arbitrary. Then
the equation~(\ref{eq0}) admits
 a two parameter group of translations only.
This group is defined by the following infinitesimal
operators:
\begin{equation} \label{op11}
X_{1}= {\ddt} , \qquad X_{2}= {\ddx} ,
\end{equation}
which generate the translations  of independent
variables.
In this case almost no limits are imposed on a mesh
and a difference equation.
In particular, we can use an orthogonal grid in the
plane $(x,t)$
which is regular in both directions, as the
conditions
(\ref{ct})--(\ref{cht})
are valid for the operators~(\ref{op11}).

The group with operators~(\ref{op11}) in the subspace
$(x,t,h,\tau,u,u_{-},u_{+},\hat{u},\hat{u}_{-},\hat{u}_{+})$
corres\-ponding to
the stencil
shown in Fig.~1
has eight  invariants:
\[
\tau, \quad h,  \quad u, \quad u_{+}, \quad u_{-},
\quad \hat{u}, \quad
\hat{u}_{-}, \quad \hat{u}_{+}.
\]
That is why any difference approximation of the
equation~(\ref{eq0}) by
the above
invariants could give difference equation which admits
the operators~(\ref{op11}).
For example, the explicit model
\begin{equation} \label{sh11}
{\hat{u}-u \over \tau } = {1 \over h}
\left( K\left( { u_{+} + u \over 2 }  \right) { \dh u{}_{x}}
- K\left( { u +  u_- \over 2 } \right) { \dh  u{}_{\bar
x} } \right)  + Q(u),
\end{equation}
where $K(u)$ and $Q(u)$ represent any approximation of
the corresponding
coefficients
by invariants and ${ { \dh u{}_{x}  } = {
u_{+} - u \over h } }$,
${ { \dh  u{}_{\bar x}  }= {u - u_{-}
\over h } }$ are
right and left difference derivatives, admits the
operators~(\ref{op11}).

{\bf 2.} If $K(u)$ is arbitrary function and $Q(u)
\equiv 0$, the equation
\begin{equation} \label{eq12}
u_{t}=(K(u)u_{x})_{x}
\end{equation}
admits a three-parameter algebra of operators (see~\cite{[18]}):
\begin{equation} \label{op12}
X_{1}={\ddt}, \qquad X_{2}= {\ddx}, \qquad X_{3}=
2t{\ddt} + x {\ddx} .
\end{equation}
This case is almost analogous to the previous one. The
operators~(\ref{op12})
do not
violate conditions of invariant orthogonality~(\ref{cht}) and
invariant uniformity of a grid (\ref{ct}), (\ref{ch}).
Thus, in this case
we could use the orthogonal grid shown in Fig.~1.
Any approximation of the equation~(\ref{eq12}) by the
seven invariants
\[
{h^{2} \over \tau }, \quad u , \quad u_{+}, \quad
u_{-}, \quad \hat{u}, \quad
\hat{u}_{-}, \quad \hat{u}_{+}
\]
gives an invariant model for the equation~(\ref{eq12}). In particular
the explicit scheme (\ref{sh11}) with $Q \equiv 0$:
\begin{equation} \label{sh12}
{\hat{u} - u \over \tau } = {1\over h}
\left( K\left( { u_{+} + u \over 2}  \right) { \dh
u{}_{x}} - K\left( { u + u_{-} \over 2 }  \right) { \dh
u{}_{\bar x} } \right)
\end{equation}
can be used.

\section{The exponential heat transfer coefficient
$\boldsymbol{K=e^u}$}

In this paragraph we consider three cases of group
classification for
$K=e^u$, in accordance
with~\cite{[4]} and~\cite{[18]}.

{\bf 1.} If $Q=0$ then the equation
\begin{equation} \label{eq21}
u_{t}=(e^u u_{x})_{x}
\end{equation}
admits a four-dimensional algebra of infinitesimal
operators:
\begin{equation} \label{op21}
X_{1}= {\ddt}, \qquad X_{2}= {\ddx}, \qquad X_{3} = 2t
{\ddt} + x{\ddx}, \qquad
X_{4}= t{\ddt} - {\ddu} .
\end{equation}
As in the cases considered above, conditions of
invariant uniformity
and invariant orthogonality are valid. A difference
model for the
equation~(\ref{eq21}) can be constructed by
approximation of
the differential equation with the help of difference
invariants:
\[
e^u { \tau \over h^{2}} ,
\quad (\hat{u} -u) , \quad (u_{+} -u), \quad (u-
u_{-}) ,
\quad (\hat{u}_{+} - \hat{u}), \quad (\hat{u} -
\hat{u}_{-}) .
\]
An example is the simple explicit difference model:
\begin{equation} \label{sh21}
{\hat{u} - u \over \tau } = {1\over h} \left( \exp
\left( {u_{+} +u \over 2}
\right) {\dh u{}_{x}} - \exp \left( { u+u_{-} \over 2}
\right) {\dh u{}_{\bar x}
} \right),
\end{equation}
but one has a lot of freedom to construct invariant
schemes
using the finite-difference invariants.

{\bf 2.} For $Q=\delta = \pm 1$ we have a possibility
to exclude the constant
source
from the equation
\begin{equation} \label{eq22}
u_{t}=(e^u u_{x})_{x} +
\delta
\end{equation}
by change of variables:
\begin{equation} \label{ch22}
\bar{u} = u - \delta t , \qquad
\bar{t} =  \delta (e^{\delta t } - 1 ) .
\end{equation}
The equation (\ref{eq22}) will be transformed into the
equation~(\ref{eq21}) by this change, but the uniformity of the
grid in the
$t$-direction is destroyed.
The equation~(\ref{eq22}) admits the four-dimensional
algebra of infinitesimal operators:
\begin{equation} \label{op22}
X_{1}= {\ddt}, \qquad X_{2}= {\ddx},
\qquad X_{3}= e^{-\delta t}{\ddt} + \delta e^{ - \delta t } {\ddu}, \qquad X_{4}=
x{\ddx} + 2{\ddu} ,
\end{equation}
and we can easily see that the operator $X_{3}$ does
not conserve
uniformity of a grid in the time direction.
The finite-difference invariants:
\[
{ e^u (e^{\delta \tau}  -1)\over h^{2} },
\quad (\hat{u} -u - \delta \tau ) , \quad (u_{+} -u),
\quad
(u- u_{-}) ,
\quad (\hat{u}_{+} - \hat{u}) , \quad (\hat{u} -
\hat{u}_{-})
\]
permit us to construct the following variant of
difference model for the
equation~(\ref{eq22}):
\begin{equation} \label{sh22}
{ \delta  (\hat{u} -u)  -  \tau \over
e^{\delta \tau} -1  } =
{1\over h} \left(
\exp \left( { u_{+} + u  \over 2 } \right)  { \dh
u{}_{x} } -
\exp \left( { u + u_{-} \over 2 } \right)  { \dh
u{}_{\bar{x}} }
\right) .
\end{equation}
Let us  note that the change of variables~(\ref{ch22})
transforms the model~(\ref{sh22}) considered on the orthogonal grid with
the time interval~$[0,T]$,
given by the formula
\begin{equation} \label{meth22}
t_{n} = \delta \ln \left(  1+  {n \over k }
\big(  e^{\delta  T } - 1 \big)
\right),\qquad n= 0,\ldots,k,
\end{equation}
where $k$ is the number of time steps of the grid,
into the model~(\ref{sh21}) with
uniform time grid on the time interval
$\big[0,\delta \left(e^{\delta T } -1 \right) \big]$.

{\bf 3.} If $Q= \pm e^{ \alpha u}$, $\alpha \neq 0$
the equation
\begin{equation} \label{eq23}
u_{t}=(e^{u} u_{x})_{x} \pm
e^{\alpha u}
\end{equation}
admits 3 infinitesimal operators:
\begin{equation} \label{op23}
X_{1}= {\ddt}, \qquad X_{2}= {\ddx},  \qquad X_{3}=
2\alpha t {\ddt} +(\alpha
-1) x{\ddx} -2{\ddu} .
\end{equation}
These operators satisfy the conditions of
orthogonality and
uniformity of invariant grids and we will consider the
stencil
of Fig.~1. Any approximation of the equation~(\ref{eq23}) by the
finite-difference invariants
\[
{\tau ^{{\alpha -1 \over 2\alpha} } \over
h },
\quad e^{\alpha u}\tau ,
\quad (\hat{u} -u), \quad (u_{+} -u), \quad (u- u_{-}),
\quad (\hat{u}_{+} - \hat{u}) , \quad (\hat{u} -
\hat{u}_{-})
\]
gives a variant of a difference model for the equation~(\ref{eq23}),
admitting the symmet\-ries (\ref{op23}),
for example, we obtain the following model:
\begin{equation} \label{sh23}
{\hat{u} -u \over  \tau} = {1\over h} \left(\exp
\left( {u_{+} +u\over 2}
\right) { \dh u{}_{x} } -
\exp \left({u+u_{-}\over2} \right) { \dh u{}_{\bar x} }
\right)  \pm    e^{\alpha u} .
\end{equation}

{\bf 4.} In accordance with group classification~\cite{[4]}
we will also consider the case $Q= \pm e^{u} + \delta$,
$\delta = \pm 1$. As in the case 2 we have the
possibility
to exclude the constant source from the equation
\begin{equation} \label{eq24}
u_{t}=(e^{u} u_{x})_{x}  \pm
 {e}^{u} +  \delta
\end{equation}
by the change of variables~(\ref{ch22}). The equation~(\ref{eq24}) will
be transformed into the equation~(\ref{eq23}). The
equation~(\ref{eq24})
admits the following infinitesimal operators:
\begin{equation} \label{op24}
X_{1}= {\ddt}, \qquad X_{2}= {\ddx}, \qquad
X_{3}= e^{ -\delta t} {\ddt} +
\delta e^{-\delta t} {\ddu} .
\end{equation}

Finite-difference invariants of (\ref{op24})
\[
e^{u} (e^{ \delta \tau} -1 ),
\quad h,
\quad (\hat{u} -u - \delta \tau ) , \quad (u_{+} -u),
\quad (u- u_{-}) ,
\quad (\hat{u}_{+} - \hat{u}) , \quad (\hat{u} -
\hat{u}_{-})
\]
permit to construct the following variant of the
difference model
\begin{equation} \label{sh24}
{ \delta  (\hat{u} -u)  -  \tau \over
e^{\delta \tau} -1  } =
{1\over h} \left(
\exp \left( { u_{+} + u  \over 2 } \right)  { \dh
u{}_{x} }-
\exp \left( { u + u_{-} \over 2 } \right)  { \dh
u{}_{\bar{x}} }
\right) \pm e^{u}
\end{equation}
on the invariant grid (\ref{meth22}).
The model for the considered equation can be
obtained from the model~(\ref{sh23}) with the help of
the
transformation~(\ref{ch22}).

\section{The power heat transfer
coefficient: $\boldsymbol{K= u^{\sigma}}$, $\boldsymbol{\sigma \neq 0, \; { - { 4\over 3}}}$}

\noindent For $K= u^{\sigma}$ further classification
depends on the source.

{\bf 1.} Let us start with the simplest case $Q \equiv
0$:
\begin{equation} \label{eq31}
u_{t}=(u^{\sigma}u_{x})_{x}.
\end{equation}

Symmetries of the equation (\ref{eq31}) are described
by the
four-dimensional algebra of infinitesimal operators
(see~\cite{[18]}):
\begin{equation} \label{op31}
X_{1}= {\ddt}, \qquad X_{2}= {\ddx} ,
\qquad  X_{3} = 2t {\ddt} + x{\ddx},
\qquad X_{4} = \sigma x{\ddx} + 2u{ \ddu} .
\end{equation}
With any $\sigma$ the operators (\ref{op31}) conserve
uniformity and
orthogonality of a grid. The finite-difference
invariants corresponding
to the stencil of Fig.~1 are
\[
u^{ \sigma } { \tau \over h^{2} },
\quad { \hat{u} \over u } ,
\quad {u_{+} \over u }, \quad {u_{-} \over u } ,
\quad {\hat{u}_{+} \over \hat{u} }, \quad {\hat{u}_{-}
\over \hat{u} } .
\]
They permit us to write, for example,
the following variant of the difference model
on the orthogonal uniform mesh in both directions:
\begin{equation} \label{sh31}
{ \hat{u} - u \over \tau } = { 1\over h }  \left(
\left( {u_{+} +u\over 2 }\right)^{\sigma} { \dh u{}_{x}}
- \left( {u +u _{-} \over 2} \right)^{\sigma} { \dh
u{}_{\bar x}}  \right) .
\end{equation}

Orthogonal grid is not the only possible way for
discrete
modeling. On the example of the equation (\ref{eq31})
we will show how to
introduce a moving mesh of the form shown in Fig.~2.
\begin{figure}[th]
{
\vspace*{-5mm}
\begin{picture}(300,160)
\put(120,0){\begin{picture}(200,150)%
{ \thicklines
\put(10,20){\vector(1,0){180}}
\put(10,20){\vector(0,1){120}}
}
\put(10,50){\line(1,0){170}}
\put(10,80){\line(1,0){160}}
\put(10,110){\line(1,0){155}}
\put(10,20){\line(1,3){10}}
\put(40,20){\line(1,3){10}}
\put(80,20){\line(1,6){5}}
\put(140,20){\line(-1,6){5}}
\put(180,20){\line(0,1){30}}
\put(20,50){\line(1,3){10}}
{\thicklines
\put(50,50){\line(1,2){15}}
\put(85,50){\line(0,1){30}}
\put(135,50){\line(-1,3){10}}
\put(50,50){\line(1,0){85}}
\put(65,80){\line(1,0){60}}
}
\put(180,50){\line(-1,3){10}}
\put(30,80){\line(1,2){15}}
\put(65,80){\line(1,4){7,5}}
\put(85,80){\line(1,6){5}}
\put(125,80){\line(-1,6){5}}
\put(170,80){\line(-1,6){5}}
\put(50,50){\circle*{5}}
\put(85,50){\circle*{5}}
\put(135,50){\circle*{5}}
\put(65,80){\circle*{5}}
\put(85,80){\circle*{5}}
\put(125,80){\circle*{5}}
\put(0,135){$t$}
\put(200,20){$x$}
\end{picture}}
%
%
\end{picture}
}
\vspace{-5mm}

\caption{A moving mesh with flat time-layers.}
\end{figure}
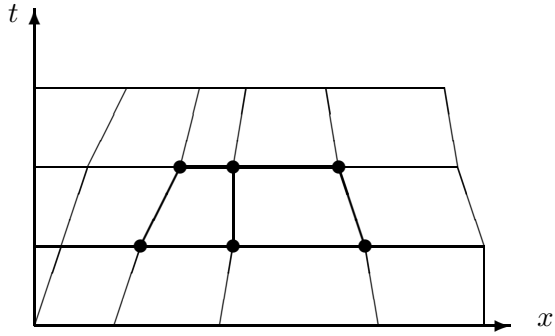

One can use an adaptive grid defined by an
evolution equation (see also~\cite{[3]})
\begin{equation} \label{ad1}
{ dx \over dt } = \varphi ( t, x, u, u_{x} ).
\end{equation}
In this case the heat transfer equation will take the
form
\begin{equation} \label{ad2}
{ du \over dt } = (u^{\sigma} u_{x} )_{x} + \varphi (
t, x, u, u_{x} ) u_{x}.
\end{equation}

Different requirements could be imposed on the
function $\varphi$.
If we require invariance of the equation~(\ref{ad1})
with respect to the
whole set of the operators~(\ref{op31}),  our freedom
to chose $\varphi$
is limited by the function
\[
\varphi = Cu^{\sigma -1} u_{x}, \qquad  C = \mbox{const}.
\]

Below we show how to introduce Lagrangian type of
evolution ${ dx \over dt }$.
Let us note that the equation (\ref{eq31}) has a
form of the conservation law that presents the
conservation of heat.
Hence we can search for a moving mesh of Lagrange type
which
evolves in accordance with heat diffusion. We should
find an evolution
${ dx \over dt}$ which satisfies the equation
\[
 { d \over  dt}
\int_{x1(t)}^{x2(t)} u\,dx = 0.
\]
Since
\[
{ d \over  dt} \int_{x1}^{x2} u\,dx  =
\int_{x1}^{x2}{ \partial u \over \partial  t}\, dx +
\left[ u {dx \over dt} \right]^{x2}_{x1} =
\left[ u^{\sigma}u_{x}  +    u {dx \over dt}
\right]^{x2}_{x1}
\]
we obtain the evolution ${dx \over dt}= - u^{\sigma
-1}u_{x}$.
Note that this evolution is invariant with respect to
the operators~(\ref{op31}). Our initial
differential equation~(\ref{eq31}) can now be
presented in the form
of the system
\begin{gather}
{dx \over dt } = - u^{\sigma -1}u_{x},\qquad
{du \over dt } =  u^{\sigma }u_{xx} +
( \sigma - 1) u^{\sigma -1} u_{x} ^2 . \label{sys31}
\end{gather}

Let us mention that the equation (\ref{eq31}) has two
conservation laws
\[
u_{t}=(u^{\sigma} u_{x} )_{x} \qquad \mbox{and} \qquad
( xu) _t =
\left( x u^{\sigma}u_{x}  - { u^{\sigma +1} \over
{\sigma +1} }  \right)_x.
\]
For the evolution system (\ref{sys31}) it is
convenient to present the conservation laws in the
integral form
\begin{equation} \label{erq1}
 { d \over  dt}
\int_{x1(t)}^{x2(t)} u\,dx = 0, \qquad
{ d \over  dt}
\int_{x1(t)}^{x2(t)} xu\,dx = -
\left.  { u^{\sigma +1} \over {\sigma +1} }  \right|
^{x_2} _{x_1}  .
\end{equation}

For difference modeling of the system (\ref{sys31}) we
can consider
the stencil shown in Fig.~3.
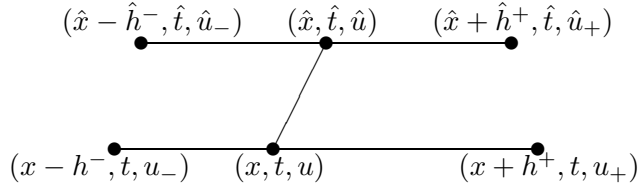
\begin{figure}[th]   \label{grnonor}
{
\vspace*{-5mm}
\begin{picture}(300,100)
\put(110,0){\begin{picture}(200,100)%
\put(20,30){\line(1,0){160}}
\put(30,70){\line(1,0){140}}
\put(80,30){\line(1,2){20}}
\put(20,30){\circle*{5}}
\put(30,70){\circle*{5}}
\put(170,70){\circle*{5}}
\put(180,30){\circle*{5}}
\put(80,30){\circle*{5}}
\put(100,70){\circle*{5}}
\put(0,75){$(\hat{x}-\hat{h}^{-},\hat{t},\hat{u}_{-})$}
\put(85,75){$(\hat{x},\hat{t},\hat{u})$}
\put(140,75){$(\hat{x}+\hat{h}^{+},\hat{t},\hat{u}_{+})$}
\put(65,20){$(x,t,u)$}
\put(150,20){$(x+h^{+},t,u_{+})$}
\end{picture}}
%
\put(90,20){$(x-h^{-},t,u_{-})$}
\end{picture}
}
\vspace{-5mm}

\caption{The stencil of the evolutionary mesh.}
\end{figure}

In the space of the variables
$(t,x,\tau,h^{+},h^{-},\hat{h}^{+},
\hat{h}^{-}, \Delta x,u,u_{+},
u_{-},\hat{u},\hat{u}_{+},\hat{u}_{-})$
correspon\-ding to this stencil there are ten
finite-difference invariants:
\[
u^{ \sigma } { \tau \over h^{+2} },
\quad { \hat{u} \over u } ,
\quad {u_{+} \over u }, \quad {u_{-} \over u } ,
\quad {\hat{u}_{+} \over \hat{u} }, \quad {\hat{u}_{-}
\over \hat{u} } ,
\quad { h^{-} \over h_{+}} ,
\quad { \hat{h}^{-} \over h^{+}} ,
\quad { \hat{h}^{+} \over h^{+}} ,
\quad { \Delta x \over h^{+}} .
\]
Approximating the system (\ref{sys31}) by these
invariants we can get, for example, the system of two
equations:
\begin{gather}
 { \Delta x \over \tau } = - { 1\over 2
\sigma }
\left(  { u_{+}^{\sigma} - u^{\sigma} \over h^{+} } +
        { u^{\sigma} - u_{-}^{\sigma} \over h^{-} }
\right) ,\qquad
 { \hat{u} + \hat{u}_+ \over 2 }
\hat{h} _+ =
  { {u} + {u}_+ \over 2 } h_+ ,\label{sh31a}
\end{gather}
where we approximated the heat conservation law to
obtain
the equation for the solution~$u$.

The first equation of system (\ref{sys31}) shows that
the evolution of~$x$ depends on the solution. The system (\ref{sh31a})
may be
inconvenient for computations
because steplength will be changed automatically and
the
nature of this process is not clear. In order to avoid
this
indeterminacy we introduce a new space variable which
values
characterize the evolution trajectories of $x$. Let us
consider
the variable $s$ defined by the system:
\[
s_{t} =  u^{ \sigma}u_{x}, \qquad
s_{x} =  u .
\]
It is easy to see that each trajectory of $x$ is
prescribed by a fixed value
of $s$ since
\[
{ds \over dt} = s_{t} + s_{x} {dx \over dt} =
u^{ \sigma}u_{x}   - u { 1\over \sigma } (
u^{\sigma})_{x} = 0.
\]
In a new coordinate system with the independent
variables $(t,s)$
the equation (\ref{eq31}) has the form
\begin{equation} \label{eq31n}
\left( { 1\over u} \right)_{t} = -
(u^{\sigma}u_{s})_{s}
\end{equation}
and the former space variable $x$ satisfies
\begin{equation} \label{sys31n}
x_{t} =  - u^{ \sigma} u_{s}, \qquad
x_{s} =  {1 \over u} .
\end{equation}
For discrete modeling of the equation (\ref{eq31}) one
can use
the equation (\ref{eq31n}) in the new independent
variables $(t,s)$
to describe the diffusion process
and the first equation of the system~(\ref{sys31n}) to
trace the evolution of the coordinate~$x$.
The equation~(\ref{eq31n}) considered together with
the system~(\ref{sys31n}) admits
the following  symmetries
\begin{gather}
X_{1}= {\ddt}, \qquad X_{2}= {\ddx} , \qquad X_{3}=
{\dds},
\qquad  X_{4} = 2t {\ddt} + s  {\dds} + x{\ddx},  \nonumber\\
X_{5} = ( \sigma + 2) s  {\dds} + \sigma x{\ddx} + 2u{
\ddu} . \label{op31n}
\end{gather}
In the new variables $(t,s)$ the stencil becomes
orthogonal
so that there
is no need to consider a nonuniform grid in the
variable~$s$.
There are following invariants for this set of
operators in the space
$(t,\tau, s , h_{s},
x,h_{x}^{+},h_{x}^{-},\hat{h}_{x}^{+},
\hat{h}_{x}^{-}, \Delta x,u,u_{+},
u_{-},\hat{u},\hat{u}_{+},\hat{u}_{-})$
corresponding to the orthogonal stencil in $(t,s)$
extended by additional
dependent variable~$x$:
\[
u^{ \sigma } { \tau \over h_{x}^{+2} },
\quad { \hat{u} \over u } ,
\quad {u_{+} \over u }, \quad {u_{-} \over u } ,
\quad {\hat{u}_{+} \over \hat{u} }, \quad {\hat{u}_{-}
\over \hat{u} } ,
\quad { h_{x}^{-} \over h_{x}^{+}} ,
\quad { \hat{h}_{x}^{-} \over h_{x}^{+}} ,
\quad { \hat{h}_{x}^{+} \over h_{x}^{+}} ,
\quad { \Delta x \over h_{x}^{+}} ,\quad { h_{s}
\over h_{x}^{+}}.
\]
By means of these invariants we get an approximation
of~(\ref{eq31n}) which has the form of a~conservation law
\begin{gather}
{ 1 \over \tau } \left( { 1 \over \hat{u} } -  { 1
\over u} \right) =
- { \alpha  \over \sigma + 1}
\left( { u_{+}^{ \sigma + 1} - 2 u^{ \sigma + 1} +
u_{-}^{ \sigma + 1}  \over h_{s}^{2} } \right)\nonumber\\
\phantom{{ 1 \over \tau } \left( { 1 \over \hat{u} } -  { 1
\over u} \right) =}{}
- { 1 -\alpha  \over \sigma + 1}
\left( { \hat{u}_{+}^{ \sigma + 1} - 2 \hat{u}^{
\sigma + 1} +
\hat{u}_{-}^{ \sigma + 1}  \over h_{s}^{2} } \right), \label{sh31n}
\end{gather}
where $0 \leq \alpha \leq 1$. Note that in the
coordinates $(t,s)$
variable~$x$ is introduced
by the system~(\ref{sys31n}) as a potential for the
equation (\ref{eq31n}).
Similarly we can introduce $x$ as a~discrete potential
with the help of the system
\begin{gather}
{ \Delta x \over \tau}
=  - { \alpha \over \sigma + 1} \left( { u_{+}^{
\sigma + 1} -
u_{-}^{ \sigma + 1}  \over 2 h_{s} } \right)
- { 1- \alpha \over \sigma + 1} \left( { \hat{u}_{+}^{
\sigma + 1} -
\hat{u}_{-}^{ \sigma + 1}  \over 2 h_{s} } \right) .\nonumber\\
{ h_{x}^{+} \over h_{s}} = {1 \over 2}
\left( {1 \over u } + { 1\over u_{+} } \right)  .
 \label{sysh31n}
\end{gather}

In computations only the equation (\ref{sh31n}) and
the
first  equation of (\ref{sysh31n}) are needed.
The second equation of the system (\ref{sysh31n})
is needed only to establish the
connection bet\-ween solutions $u(x)$ and $u(s)$ for
a fixed time.
Given some initial data  $ u ( 0, x) =   u_0 ( x)$,
we choose an appropriate steplength $h_s$ for the
Lagrangian mass
coordinate $s$. Then we can introduce the mesh points
$x_i$ in the
original coordinates satisfying
\[
{ x_{i+1} - x_{i} \over h_s } =
{1 \over 2 }
\left(
{ 1 \over u_0 (x_i) } +  { 1 \over u_0 (x_{i+1}) }
\right),
\]
i.e., we use this equation to establish difference
relation between
the original space coordinate~$x$ and the Lagrangian
mass
coordinate~$s$. Computing the solution with the help
of the
numerical scheme~(\ref{sh31n}) and the first equation
of~(\ref{sysh31n}),
we preserve the relation
\[
{ x_{i+1} - x_{i} \over h_s } =
{1 \over 2 } \left( { 1 \over u_i } +  { 1 \over
u_{i+1} } \right).
\]

Introducing the mass type variable $s$,  we
can rewrite the conservation laws (\ref{erq1}) as
\[
 { \partial  \over  \partial t}
\int_{s1}^{s2} ds = 0 , \qquad
{ \partial  \over  \partial t}
\int_{s1}^{s2} x\,ds = -
\left.  { u^{\sigma +1} \over {\sigma +1} }  \right|
^{s_2} _{s_1}  .
\]
The proposed discrete model possesses difference
analogs of
these conservation laws
\begin{gather*}
\sum_{i =1}^{N-1}  h_s = {\rm const}, \\
\sum_{i =1}^{N-1} { \hat{x}_{i} + \hat{x}_{i+1} \over
2}   \,h_s
- \sum_{i =1}^{N-1} { {x}_{i} + {x}_{i+1} \over 2}\,
h_s =
- { \alpha \over \sigma + 1}
\left( { u_{N+1}^{ \sigma + 1} + u_{N}^{ \sigma + 1}
\over 2 } \right)
\\
\qquad{}- {  1 - \alpha \over \sigma + 1}
\left( { \hat{u}_{N+1}^{ \sigma + 1} +
\hat{u}_{N}^{ \sigma + 1}  \over 2  } \right)
+
 { \alpha \over \sigma + 1} \left( { u_{-1}^{ \sigma +
1} +
u_{0}^{ \sigma + 1}  \over 2 } \right)
+ {  1 - \alpha \over \sigma + 1} \left( {
\hat{u}_{-1}^{ \sigma + 1} +
\hat{u}_{0}^{ \sigma + 1}  \over 2  } \right).
\end{gather*}
Let us mention that for computations we need to
propose some
method for the boundary points.

{\bf 2.} $Q=\delta u$,  $ \delta = \pm 1$. In this
case the
symmetry of the equation
\begin{equation} \label{eq32}
u_{t}=(u^{\sigma}u_{x})_{x} + \delta u
\end{equation}
is described by the following infinitesimal operators:
\begin{gather}
X_{1}= {\ddt}, \qquad X_{2}= {\ddx}, \qquad X_{3}=
\sigma x {\ddx} + 2u{\ddu},
\nonumber\\
X_{4}= e^{- \delta \sigma t}
{\ddt} + \delta e^{- \delta \sigma t}
u{\ddu} .\label{op32}
\end{gather}
By the change of variables
\begin{equation} \label{ch32}
\bar{u} = u e^{ - \delta t }, \qquad
\bar{t} = { \delta \over \sigma }
(e^{\delta \sigma t } - 1 )
\end{equation}
the equation (\ref{eq32}) is transformed into the
equation (\ref{eq31}).

The finite-difference invariants
\[
{ u^{\sigma}  ( e^{ \delta \sigma \tau} -1 )
\over h^{2}},
\quad \left( \delta \ln {\hat{u} \over u} -  \tau
\right) ,
\quad {u_{+} \over u},
\quad {u_{-} \over u} ,
\quad {\hat{u}_{+} \over  \hat{u}} ,
\quad {\hat{u}_{-} \over \hat{u} }
\]
give the following possibility for an explicit
difference model
\begin{equation} \label{sh32}
{ \sigma u \over
e^{ \delta \sigma \tau} -1  }
\left( \delta \ln {\hat{u} \over u} -  \tau \right) =
{1\over h} \left(
\left( {u_{+} +  u \over 2} \right)^{\sigma} { \dh
u{}_{x}} -
\left( { u + u_{-} \over 2} \right)^{\sigma} { \dh
u{}_{\bar{x}} }\right) .
\end{equation}
Let us  remark that the change of variables
(\ref{ch32}) transforms
this equation considered on the orthogonal  mesh with
time layers
\begin{equation} \label{meth32}
t_{n} = { \delta \over \sigma } \ln \left( 1 + {n
\over k} \big(
e^{\delta \sigma T} -1\big)
\right),\qquad
  n = 0,\ldots,k,
\end{equation}
which fill the time interval $[0,T]$,
into the equation~(\ref{sh31}) on
the uniform grid on the time interval
$\left[0, { \delta \over \sigma }
\big( e^{ \delta \sigma T} -1\big)\right]$.

{\bf 3.} $Q= \pm u^{\sigma +1} + \delta u ^n$,  $
\delta = \pm 1$.
The equation
\begin{equation} \label{eq33}
u_{t}=(u^{\sigma}u_{x})_{x}  \pm
u^{n},\qquad\sigma,n={\rm const},
\end{equation}
admits a three-parameter symmetry group. A possible
representation
of this group is by the
following infinitesimal operators:
\begin{equation} \label{op33}
X_{1}= {\ddt}, \qquad X_{2}={\ddx}, \qquad X_{3}=
2(n-1)t{\ddt}+(n-\sigma
-1)x{\ddx} -2u{\ddu} .
\end{equation}
The set (\ref{op33}) satisfies all the conditions
(\ref{ct})--(\ref{cht}).
So we can use an orthogonal grid that is uniform in
both~$t$ and~$x$
directions.
By considering the set of the operators (\ref{op33})
in the space
$( t, \hat{t}, x, h^{+}, h^{-}, u, u_{+},
u_{-}, \hat{u}, \hat{u}_{+}, \hat{u}_{-})$ that
corresponds to the stencil shown in Fig.~1 we find~7 difference invariants of the Lie algebra:
\[
{\tau ^{ n -\sigma -1 \over 2(n-1) } \over h },
\quad \tau u^{n-1} , \quad { \hat{u} \over u} ,
\quad {u_{+} \over u}, \quad { u_{-} \over u} ,
\quad {\hat{u}_{+} \over \hat{u}}, \quad { \hat{u}_{-}
\over \hat{u}} .
\]
The small number of symmetry operators provides us
with a large
number of
difference invariants. Thus we are left with some
additional
degrees of freedom in invariant difference
modeling of~(\ref{eq33}). By means of the discrete
invariants we obtain
the following explicit scheme:
\begin{equation} \label{sh33}
{ \hat{u} - u \over \tau } =
{ 1  \over h} \left(  \left({u_{+} +u \over 2}
\right)^{\sigma} { \dh  u{}_{x} }
- \left({u+ u_{-}\over 2 }\right)^{\sigma} {\dh
u{}_{\bar x}}  \right)  \pm
  u^{n},
\end{equation}
where $\dh u{}_{x}  = {u_{+} - u \over h
}$, $\dh u{}_{ \bar x}  =
{u - u_{-} \over h}$.

{\bf 4.} $Q= \pm u^{\sigma +1} + \delta u$,  $ \delta
= \pm 1$. The equation
\begin{equation} \label{eq34}
u_{t}=(u^{\sigma}u_{x})_{x} \pm  u^{\sigma+1} + \delta
u
\end{equation}
is connected with the equation (\ref{eq33}) by the
transformation (\ref{ch32}).
The infinitesimal operators admitted by the equation
are
\begin{equation} \label{op34}
X_{1}= {\ddt}, \qquad X_{2}= {\ddx},
\qquad X_{3} = e^{-  \delta
\sigma  t} {\ddt} + {\delta} e^{- \delta  \sigma t} u {\ddu} .
\end{equation}
With the help of invariants for the operators
(\ref{op34})
\[
u^{\sigma}  ( e^{ \delta
\sigma \tau} -1 ),
\quad h,
\quad \left( \delta \ln {\hat{u} \over u} -  \tau
\right) ,
\quad {u_{+} \over u},
\quad {u_{-} \over u} ,
\quad {\hat{u}_{+} \over  \hat{u}} ,
\quad {\hat{u}_{-} \over \hat{u} }
\]
we have the following example of an explicit
difference model:
\begin{equation} \label{sh34}
{ \sigma u \over
e^{ \delta \sigma \tau} -1  }
\left( \delta \ln {\hat{u} \over u} -  \tau \right) =
{1\over h} \left(
\left( { u_{+} +  u \over 2} \right)^ {\sigma } { \dh
u{}_{x}} -
\left({ u + u_{-} \over 2} \right) ^{\sigma } { \dh
u{}_{\bar{x}} } \right)
\pm u^{\sigma+1}.
\end{equation}
This equation considered on the grid (\ref{meth32})
is transformed into the equation~(\ref{sh33}) on a~uniform time grid by the variable
change~(\ref{ch32}).

\section{The special case of power
heat transfer\\
coefficient: $\boldsymbol{K= u^{ - 4/3}}$}

{\bf  1.} If $Q \equiv 0$, then the symmetry of the
equation
\begin{equation} \label{eq41}
u_{t}=\big(u^{- 4/3}u_{x}\big)_{x}
\end{equation}
is described by the five-dimensional algebra of
infinitesimal operators (see~\cite{[18]}):
\begin{gather}
 X_{1}= {\ddt}, \qquad X_{2}= {\ddx} ,
\qquad  X_{3} = 2t {\ddt} + x{\ddx}, \nonumber \\
X_{4} = 2 x{\ddx} - 3u{ \ddu} ,
\qquad
X_{5}= x^{2} {\ddx} -3xu{\ddu} .
\label{op41}
\end{gather}
These operators conserve orthogonality and uniformity
of
a grid in the time direction. The operator $X_{5}$
conserve
uniformity in the $t$-direction, but does not conserve
uniformity of the grid
in the $x$-direction;
however orthogonality is not disturbed. We will
consider the
stencil shown in Fig.~4.
\begin{figure}[th]
{
\vspace*{-5mm}
\begin{picture}(300,100)
\put(100,0){\begin{picture}(200,100)%
\put(30,30){\line(1,0){160}}
\put(30,70){\line(1,0){160}}
\put(90,30){\line(0,1){40}}
\put(30,30){\circle*{5}}
\put(30,70){\circle*{5}}
\put(190,70){\circle*{5}}
\put(190,30){\circle*{5}}
\put(90,30){\circle*{5}}
\put(90,70){\circle*{5}}
\put(74,75){$(x,\hat{t},\hat{u})$}
\put(152,75){$(x+h^{+},\hat{t},\hat{u}_{+})$}
\put(74,20){$(x,t,u)$}
\put(152,20){$(x+h^{+},t,u_{+})$}
\end{picture}}
\put(92,75){$(x-h^{-},\hat{t},\hat{u}_{-})$}
\put(92,20){$(x-h^{-},t,u_{-})$}
\end{picture}
}
\vspace{-5mm}
\caption{The stencil of nonuniform mesh.}
\end{figure}
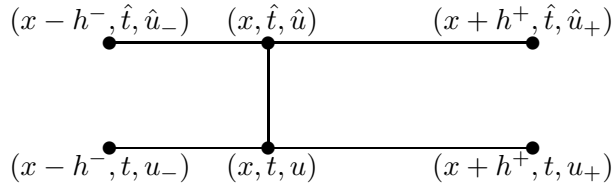

 The finite-difference
invariants
corresponding to this stencil
\[
{\hat{u} \over u}, \quad {\hat{u}_{+} \over u_{+}},
\quad {\hat{u}_{-} \over u_{-}},
\quad u_{+}^{1/3}u^{1/3} { h^{+} \over \sqrt{\tau} },
\quad u_{-}^{1/3}u^{1/3} { h^{-} \over \sqrt{\tau} },
\quad { u^{2/3} \over \sqrt{\tau}}
\left( { h^{+}  h^{-} \over h^{+} + h^{-} } \right)
\]
give among others the explicit difference model
\begin{equation} \label{sh41}
{\hat{u} - u \over \tau } =
-   { h^{+} + h^{-} \over 6 h^{+} h^{-} }
\left(
{ u_{+}^{-1/3}  - u^{-1/3} \over h^{+} } -
{ u^{-1/3}  - u_{-}^{-1/3} \over h^{-} }
\right).
\end{equation}

\begin{remark} Let us note that we can
not propose a space
mesh $ h_{+} = f ( x , h_{-} )$ which is preserved
under
all transformations of the group (\ref{op41}). It can
be clearly seen
from the absence of difference invariants in the space
$ ( x, h_{-}, h_{+} )$. For example, if we take a
solution of the
difference scheme (\ref{sh41}) on a regular mesh
$ h_{-} = h_{+} $, a general group transformation
corresponding to~(\ref{op41}) will transform the solution into another
solution of
this  difference scheme but possibly on a nonuniform
mesh. This remark
is also valid for the cases~{2},~{4} and~{5} of this
section.
\end{remark}

{\bf 2.} In the case $Q=\delta u$,  $ \delta = \pm 1$
equation
\begin{equation} \label{eq42}
u_{t}=\big(u^{-4/3}u_{x}\big)_{x} + \delta u
\end{equation}
admits operators
\begin{gather}
X_{1}= {\ddt}, \qquad X_{2}= {\ddx},
\qquad X_{3}= 2x {\ddx} -
3u {\ddu},  \nonumber\\
X_{4}= e^{{ 4 \delta  t \over 3}}
{\ddt} + \delta e^{ { 4 \delta  t \over 3 }} u{\ddu},
\qquad
X_{5}= x^{2} {\ddx} - 3xu {\ddu} . \label{op42}
\end{gather}

By the change of variables (\ref{ch32}) this equation
can be transformed
into the equation~(\ref{eq41}).
Let us write out the difference invariants
for the set of the operators~(\ref{op42}):
\begin{gather*}
\left( \delta \ln {\hat{u} \over u} -  \tau
\right) ,
\quad u^{2/3} \left( { h^{+}  h^{-} \over h^{+} +
h^{-} } \right)
{1\over \sqrt{
( e^{ \delta \sigma \tau} -1 )
}}   , \\
{ u_{+}^{1/3}u^{1/3}  h^{+} \over
{ \sqrt{
(e^{ \delta \sigma \tau} -1 )
}}  },
\quad { u_{-}^{1/3}u^{1/3}  h^{-} \over
{ \sqrt{
( e^{ \delta \sigma \tau} -1 )
}}  },
\quad  { \hat{u}_{+}^{1/3}\hat{u}^{1/3}  h^{+} \over
{ \sqrt{
( e^{  \delta \sigma \tau} -1 )
}}  },
\quad { \hat{u}_{-}^{1/3}\hat{u}^{1/3}  h^{-} \over
{ \sqrt{
( e^{\delta \sigma \tau} -1 )
}}  }.
\end{gather*}
These invariants can be used for construction of a
difference
model for the equation (\ref{eq42}). We show the
explicit variant
of the difference model:
\begin{equation} \label{sh42}
{ \sigma u \over
e^{\delta \sigma \tau} -1  }
\left( \delta \ln {\hat{u} \over u} -  \tau \right) =
-  { h^{+} + h^{-} \over 6 h^{+} h^{-} }
\left(
{ u_{+}^{-1/3}  - u ^{-1/3} \over h^{+} } -
{ u^{-1/3}  - u_{-}^{-1/3} \over h^{-} }
\right).
\end{equation}

{\bf 3.} $Q= \pm u^{n}$,  $n \neq { -{1\over 3}}$. The
equation
\begin{equation} \label{eq43}
u_{t}=\big(u^{-  4/3 } u_{x}\big)_{x} \pm  u^{n}
\end{equation}
admits infinitesimal operators:
\begin{equation} \label{op43}
X_{1}= {\ddt}, \qquad X_{2}={\ddx}, \qquad X_{3}=
2(n-1)t{\ddt}+(n + {1\over 3}
)x{\ddx} -2u{\ddu} .
\end{equation}
Although this equation is specified in group
classification (see~\cite{[4]}),
it is a particular case of the equation~(\ref{eq33}),
--- there is no
extension of the admitted group. That's why as an
invariant difference model
for the equation (\ref{eq43}) we can use the model~(\ref{sh33})
with parameter $\sigma = - { 4\over 3}$,
corresponding to the given equation.

{\bf 4.} If $Q= \alpha  u^{- 1/ 3}$,  $\alpha =
\pm 1$, then the
variant
of the difference model for the equation
\begin{equation} \label{eq44}
u_{t}=\big(u^{-  4/3 } u_{x}\big)_{x} \pm  u^{- 1/3 }
\end{equation}
depends on the sign of the coefficient $\alpha$. The
equation~(\ref{eq44})
admits a five-dimensional algebra of infinitesimal
operators, namely
\begin{gather}
X_{1}= {\ddt}, \qquad X_{2}={\ddx}, \qquad X_{3}=
{4\over 3} t{\ddt} + 2u{\ddu}
,\nonumber \\
X_{4}= e^{2 \sqrt{ \alpha / 3
} x } {\ddx} -
\sqrt{3 \alpha } e^{2 \sqrt{
\alpha / 3 } x } u {\ddu}
,   \nonumber\\
X_{5}= e^{ - 2 \sqrt{ \alpha /
3 } x } {\ddx} +
\sqrt{3 \alpha } e^{ - 2
\sqrt{ \alpha / 3 } x } u
{\ddu}  .
\label{op44}
\end{gather}
a.) The case $\alpha = 1$.
By the change of variables
\begin{equation} \label{ch44a}
\bar{u} = u \cosh^{3}
\left(  { x \over \sqrt{3} } \right), \qquad
\bar{x} =   \sqrt{  3  }
\tanh \left(  { x \over \sqrt{3} }  \right)
\end{equation}
we transfer the considered equation into the equation~(\ref{eq41})
(see~\cite{[1]}).
With the help of difference
invariants
\begin{gather*}
{ \hat{u} \over u} ,
\quad { \hat{u}_{+} \over u_{+} } ,
\quad { \hat{u}_{-} \over u_{-} } ,
\quad
 \sqrt{\tau}  u^{-2/3}
\left(
 { 1 \over \tanh \left( { h^{+}
\over  \sqrt{3} }  \right) } +
 { 1 \over \tanh \left({ h^{-}
\over  \sqrt{3} }  \right) }
\right),\\
{ 1 \over \sqrt{\tau} } u^{1/3} u_{+}^{1/3}
\sinh \left( { h^{+} \over  \sqrt{3} } \right) ,
\quad
{ 1 \over \sqrt{\tau} } u^{1/3} u_{-}^{1/3}
\sinh \left( { h^{-} \over  \sqrt{3} } \right)
\end{gather*}
we can construct a difference model. Let us write out
one of the possible
variants of the difference model, namely an explicit
model:
\begin{gather}
{\hat{u} - u \over \tau } =
 - { 1 \over 18}
\left(
 { 1 \over \tanh
\left(  { h^{+} \over  \sqrt{3} }
\right) } +
 { 1 \over \tanh
\left(  { h^{-} \over  \sqrt{3} }
\right) }
\right)
\nonumber\\
\phantom{{\hat{u} - u \over \tau } =}{}\times
 \left(
{  u_{+}^{-1/3} - u^{-1/3} \cosh
\left( { h^{+} \over  \sqrt{3} }
 \right)
\over \sinh
\left(  { h^{+} \over  \sqrt{3} } \right) } -
{ u^{-1/3} \cosh
\left(  { h^{-} \over  \sqrt{3} }
\right) -
u_{-}^{-1/3}
\over \sinh
\left(  { h^{-} \over  \sqrt{3} }
\right) }
\right).
\label{sh44a}
\end{gather}

\noindent b.) The case $\alpha = - 1$.
By the change of variables
\begin{equation} \label{ch44b}
\bar{u} = u \cos^{3} \left(  { x \over
\sqrt{3} }  \right), \qquad
\bar{x} =   \sqrt{  3  }
 \tan \left(  { x \over \sqrt{3}
}   \right)
\end{equation}
we can transfer the given equation into the equation~(\ref{eq41})
(see~\cite{[1]}).
The set of finite-difference invariants:
\begin{gather*}
{ \hat{u} \over u} ,
\quad { \hat{u}_{+} \over u_{+} } ,
\quad { \hat{u}_{-} \over u_{-} } ,
\quad
 \sqrt{\tau}  u^{-2/3}
\left(
 { 1 \over \tan
\left( { h^{+} \over  \sqrt{3} }
\right) } +
 { 1 \over \tan
\left(  { h^{-} \over  \sqrt{3} } \right) }
\right),\\
{ 1 \over \sqrt{\tau} } u^{1/3} u_{+}^{1/3}
\sin \left( { h^{+} \over  \sqrt{3} } \right) ,
\quad
{ 1 \over \sqrt{\tau} } u^{1/3} u_{-}^{1/3}
\sin \left( { h^{-} \over  \sqrt{3} } \right)
\end{gather*}
provides us with a possibility to construct an
invariant difference
scheme. For example, one can use an explicit
difference model:
\begin{gather}
{\hat{u} - u \over \tau } =
 - { 1 \over 18}
\left( { 1 \over \tan \left( { h^{+} \over  \sqrt{3} }
\right) } +  { 1 \over \tan \left(  { h^{-} \over  \sqrt{3} }
\right) } \right) \nonumber\\
\phantom{{\hat{u} - u \over \tau } =}{}
\times  \left({  u_{+}^{-1/3} - u^{-1/3} \cos \left(
{ h^{+} \over  \sqrt{3} }  \right)
\over \sin \left( { h^{+} \over \sqrt{3} }  \right) } -
{ u^{-1/3} \cos \left(  { h^{-} \over \sqrt{3} }  \right) -
u_{-}^{-1/3} \over \sin \left( { h^{-} \over \sqrt{3} }  \right) } \right).\label{sh44b}
\end{gather}
We stress that the obtained difference models
(\ref{sh44a}) and~(\ref{sh44b})
are connected with the difference model~(\ref{sh41})
for the equation~(\ref{eq41})
by the changes of variables~(\ref{ch44a}) and~(\ref{ch44b}) correspondingly
as the initial differential equations.

{\bf 5.} $Q= \alpha  u^{ - 1/3} + \delta u$,
$|\alpha| = |\delta| =
1$.
As in the previous point, two cases of parameter
$\alpha$ in the equation
\begin{equation} \label{eq45}
u_{t}=(u^{\sigma}u_{x})_{x} \pm  u^{\sigma+1} + \delta
u
\end{equation}
should be considered separately and two difference
models should be
constructed.
Let us write out the infinitesimal operators,
admitted by the equation~(\ref{eq45}):
\begin{gather}
X_{1}= {\ddt}, \qquad X_{2}= {\ddx},
\qquad X_{3} = e^{ 4  \delta
t\over 3 } {\ddt} +
{\delta} e^{ {4 \delta t\over 3}}
u {\ddu} ,\nonumber\\
X_{4}= e^{ 2 \sqrt{ \alpha / 3
} x } {\ddx} -
\sqrt{3 \alpha } e^{2 \sqrt{
\alpha / 3 } x } u {\ddu},\nonumber\\
X_{5}= e^{ - 2 \sqrt{ \alpha /
3 } x } {\ddx} +
\sqrt{3 \alpha } e^{   - 2
\sqrt{ \alpha / 3 } x } u
{\ddu} . \label{op45}
\end{gather}
\noindent  a.) The case $\alpha = 1$.
The change of variables (\ref{ch32}) transforms
the considered equation into the equation~(\ref{eq44}) 
and the change (\ref{ch44a}) into the equation~(\ref{eq41}).

We write out the set of finite-difference invariants
for the equation~(\ref{eq45})
with $\alpha = 1$:
\begin{gather*}
\left( \delta \ln {\hat{u} \over u} -  \tau \right) ,
\quad
\sqrt{
( e^{\delta \sigma \tau} -1 )
}
u^{-2/3}
\left(
 { 1 \over \tanh \left(  { h^{+}
\over  \sqrt{3} }  \right) } +
 { 1 \over \tanh \left(  { h^{-}
\over  \sqrt{3} }  \right) }
\right),\\
{ u^{1/3} u_{+}^{1/3}
\over
\sqrt{ ( e^{  \delta \sigma
\tau} -1 ) } }
\sinh \left( { h^{+} \over  \sqrt{3} } \right) ,
\quad
{ u^{1/3} u_{-}^{1/3}
\over
\sqrt{ ( e^{ \delta \sigma
\tau} -1 ) } }
\sinh \left( { h^{-} \over  \sqrt{3} } \right) ,
\\
{ \hat{u}^{1/3} \hat{u}_{+}^{1/3}
\over
\sqrt{ ( e^{ \delta \sigma
\tau} -1 ) } }
\sinh \left( { h^{+} \over  \sqrt{3} } \right) ,
\quad
{ \hat{u}^{1/3} \hat{u}_{-}^{1/3}
\over
\sqrt{ ( e^{ \delta \sigma
\tau} -1 ) } }
\sinh \left( { h^{-} \over  \sqrt{3} } \right) .
\end{gather*}
The explicit variant of the difference model for the
equation (\ref{eq45}) on
the time grid (\ref{meth32}) has the form:
\begin{gather}
{ \sigma u \over
e^{  \delta \sigma \tau} -1  }
\left( \delta \ln {\hat{u} \over u} -  \tau \right) =
 - { 1 \over 18}\left(
 { 1 \over \tanh \left({ h^{+}
\over  \sqrt{3} }  \right) } +
 { 1 \over \tanh \left({ h^{-}
\over  \sqrt{3} }  \right) }
\right) \nonumber\\
\qquad{}\times \left(
{  u_{+}^{-1/3} - u^{-1/3} \cosh
\left(  { h^{+} \over  \sqrt{3} }
 \right)
\over \sinh
\left( { h^{+} \over  \sqrt{3} }
\right) } -
{ u^{-1/3} \cosh
\left( { h^{-} \over  \sqrt{3} }
\right) -
u_{-}^{-1/3}
\over \sinh
\left(  { h^{-} \over  \sqrt{3} }
\right) }
\right). \label{sh45a}
\end{gather}

\noindent b.) The case $\alpha = - 1$
Using the change of variables (\ref{ch44b})
we can transfer this equation into the equation~(\ref{eq42}) 
and by the change~(\ref{ch32}) into the equation~(\ref{eq44}).
Difference model for the equation~(\ref{eq45}) can be
obtained with
the help of the invariants
\begin{gather*}
\left( \delta \ln {\hat{u} \over u} -  \tau \right) ,
\quad
\sqrt{
( e^{ \delta \sigma \tau} -1 )
}
u^{-2/3}
\left(
 { 1 \over \tan
\left(  { h^{+} \over  \sqrt{3} }
\right) } +
 { 1 \over \tan
\left(  { h^{-} \over  \sqrt{3} }
\right) }
\right),
\\
{ u^{1/3} u_{+}^{1/3}
\over
\sqrt{ (e^{  \delta \sigma
\tau} -1 ) } }
\sin \left( { h^{+} \over  \sqrt{3} } \right) ,
\quad
{ u^{1/3} u_{-}^{1/3}
\over
\sqrt{ ( e^{ \delta \sigma
\tau} -1 ) } }
\sin \left( { h^{-} \over  \sqrt{3} } \right) ,
\\
{ \hat{u}^{1/3} \hat{u}_{+}^{1/3}
\over
\sqrt{ (e^{ \delta \sigma
\tau} -1 ) } }
\sin \left( { h^{+} \over  \sqrt{3} } \right) ,
\quad
{ \hat{u}^{1/3} \hat{u}_{-}^{1/3}
\over
\sqrt{ ( e^{  \delta \sigma
\tau} -1 ) } }
\sin \left( { h^{-} \over  \sqrt{3} } \right) .
\end{gather*}
One of possible difference models for the equation~(\ref{eq45}) is
\begin{gather}
{ \sigma u \over
e^{ \delta \sigma \tau} -1  }
\left( \delta \ln {\hat{u} \over u} -  \tau \right) =
 - { 1 \over 18}
\left(
 { 1 \over \tan
\left( { h^{+} \over  \sqrt{3} }
\right) } +
 { 1 \over \tan
\left(  { h^{-} \over  \sqrt{3} }
\right) }
\right)\nonumber\\
\qquad{} \times
 \left(
{  u_{+}^{-1/3} - u^{-1/3} \cos \left(
{ h^{+} \over  \sqrt{3}}  \right)
\over \sin \left( { h^{+} \over
\sqrt{3} }  \right) } -
{ u^{-1/3} \cos \left(  { h^{-} \over
\sqrt{3} }  \right) -
u_{-}^{-1/3}
\over \sin \left(  { h^{-} \over
\sqrt{3} }  \right) }.
\right). \label{sh45b}
\end{gather}
The difference
models (\ref{sh45a}) and (\ref{sh45b}) are connected
with the model~(\ref{sh42})
by changes of variables (\ref{ch44a}) and~(\ref{ch44b}) correspondingly.
The variable change~(\ref{ch32}) transforms the
difference models obtained in
this point
into model of the point~4 for corresponding values
of the parameter $\alpha$.
This example shows that
in invariant difference modeling it is possible to get
consistent models
which are connected with each other by the same point
transformations
as their original differential counterparts.

\section{Linear heat conductivity with a source}

In this section we consider the semilinear heat
transfer equation
\begin{equation}  \label{eq5}
u_{t}=u_{xx} + Q(u)
\end{equation}
with different types of a source.

{\bf 1.} With $Q= \pm e^{u}$
the equation becomes
\begin{equation} \label{eq51}
u_{t}=u_{xx} \pm e^{u}.
\end{equation}
It admits a three-dimensional algebra of
infinitesimal operators
\begin{equation} \label{op51}
X_{1}= {\ddt}, \qquad X_{2}= {\ddx}, \qquad X_{3}=
2t{\ddt}+x{\ddx}-2{\ddu} .
\end{equation}
It is easy to check that for the operators the
conditions of
orthogonality and uniformity conservation of a grid
hold.
Approximation of the equation by the invariants
\[
{ h^{2} \over \tau }, \quad \tau e^{u},
 \quad (\hat {u} -u) , \quad (u_{+} -u) , \quad (  u-
u_{-}) ,
 \quad (\hat {u}_{+} - \hat{u}) , \quad (\hat{u} -
\hat{u}_{-})
\]
will give different types of difference models. An
explicit one is
\begin{equation} \label{sh51}
{\hat{u} - u \over \tau} = {1 \over h } ( { \dh u{}_{x}}
- { \dh u{}_{\bar x}} )
\pm e^{u}.
\end{equation}

{\bf 2.} $Q=\pm u^{n}$. The equation
\begin{equation} \label{eq52}
u_{t}=u_{xx} \pm  u^{n}
\end{equation}
admits the following infinitesimal operators:
\begin{equation} \label{op52}
X_{1}= {\ddt}, \qquad X_{2}= {\ddx}, \qquad X_{3}=
2(n-1)t{\ddt} + (n-1) x{\ddx}
- 2u{\ddu} .
\end{equation}
The operators satisfy to conditions
(\ref{ct})--(\ref{cht}) and the
finite-difference invariants
\[
{ h^{2} \over  \tau } , \quad \tau u ^{n-1},
\quad {\hat {u} \over u}, \quad {u_{+} \over u}, \quad
{ u_{-} \over u} ,
\quad {\hat {u}_{+} \over \hat{u} }, \quad {
\hat{u}_{-} \over \hat{u}}
\]
permit us construct, for example,
the following difference scheme:
\begin{equation} \label{sh52}
{\hat{u} - u \over \tau } = {1 \over h }  ( { \dh
u{}_{x}} - { \dh  u{}_{\bar x}} )
 \pm u^{n}.
\end{equation}

{\bf 3.} $Q= \delta u \ln u$, $\delta = \pm 1$.
 The semilinear heat transfer
equation
\begin{equation} \label{eq53}
u_{t}=u_{xx} + \delta  u \ln u, \qquad  \delta = \pm 1
\end{equation}
admits the four-parameter Lie symmetry group of point
transformations~\cite{[4]} corresponding to the
following
set of infinitesimal operators:
\begin{equation} \label{op53}
X_{1}= {\ddt}, \qquad X_{2} = {\ddx}, \qquad  X_{3} = 2
e^{ \delta t } {\ddx} - \delta e^{\delta t} x u{\ddu},
\qquad X_{4} = e^{\delta t } u
{\ddu} .
\end{equation}
Before constructing a difference equation and a grid
that
approximate~(\ref{eq53}) and inherit the whole Lie
algebra~(\ref{op53})
we should first check condition~(\ref{cht}) for
the invariance of orthogonality. The operators
$X_{1}$, $X_{2}$
and $X_{4}$ conserve orthogonality, while $X_{3}$ does
not: the condition
(\ref{cht}) is not true for operator $X_{3}$.
Consequently an orthogonal
mesh can not be used for the invariant modeling of~(\ref{eq53}).
Conditions (\ref{cht2}) are true for the complete set
of operators, so
it is possible to use a nonorthogonal grid with flat
time layers
and we will use the grid shown in Fig.~2.

A possible reformulation of equation (\ref{eq53}) by
using the four
differential invariants in the subspace
$(t, x, u, u_{x}, u_{xx}, dt, dx, du)$:
\begin{gather*}
J_{1} = dt, \qquad
J_{2} = \left( { u_{x} \over u } \right) ^{2} - {
u_{xx}
\over u } ,
\qquad
J_{3} = 2 { u_{x} \over u } + { dx \over dt } ,
\\
J_{4} = {du \over u dt } - \delta  \ln u + { 1\over 4}
\left( {dx \over dt}
\right)^{2}
\end{gather*}
is given by the system
\[
J_{3} = 0, \qquad J_{4} = J_{2}
\]
that is
\begin{equation} \label{sys53}
{dx\over dt}  = -2 {u_{x} \over u }, \qquad
{du\over dt}  = u_{xx} + \delta u \ln u
- 2 { u_{x}^{2} \over u
}.
\end{equation}
So, the structure of the admitted group suggests the
use of two
evolution equations.

As the next step, we will find difference invariants
for the set
$X_{1}$--$X_{4}$ of the group~(\ref{op53}). These
invariants are necessary
for the approximation of system (\ref{sys53}). We will
use
the six-point difference stencil of Fig.~3
on which we will approximate system (\ref{sys53}). The
stencil
defines the difference subspace $( t, \hat{t}, x,
\hat{x}, h^{+}, h^{-},
\hat{h}^{+},
\hat{h}^{-}, u, u_{+}, u_{-},\hat{u},
\hat{u}_{+}, \hat{u}_{-} )$
and the group~(\ref{op53}) has the following
difference invariants
\begin{gather*}
I_{1} =  \tau, \quad I_{2} =  {h}^{+}, \quad I_{3} =
{h}^{-},
\quad I_{4} = \hat{h}^{+}, \quad I_{5} = \hat{h}^{-},
\\
I_{6} = ( \ln u) _{x}  - ( \ln u )_{\bar{x}} ,
\quad I_{7} = ( \ln \hat{u} ) _{x}  - ( \ln \hat{u}
)_{\bar{x}} ,
\\
I_{8} = \delta \Delta x + 2 (e^{\delta \tau} - 1)
\left( { h^{-} \over h^{+} + h^{-} } ( \ln u )_{x}  +
    { h^{+} \over h^{+} + h^{-} } (  \ln u )_{\bar{x}}
 \right) ,
\\
I_{9} = \delta \Delta x + 2 (1 - e^{- \delta \tau} )
\left( { \hat{h}^{-} \over \hat{h}^{+} + \hat{h}^{-} }
( \ln \hat{u} )_{x} +
   {\hat{h}^{+} \over \hat{h}^{+} + \hat{h}^{-} } (
\ln \hat{u} )_{\bar{x}}
\right) ,
\\
I_{10} = \delta (\Delta x) ^{2} + 4 (1 - e^{- \delta
\tau} )
\left( \ln \hat{u} - e^{ \delta
\tau} \ln u \right) ,
\end{gather*}
where $\Delta x = \hat{x}-x $,
$ (\ln u)_{x}={ \ln u_{+} - \ln u \over
h^{+}}$ ,
$( \ln u ) _{\bar x} = { \ln u - \ln
u_{-} \over h^{-}} $.

 An explicit model can be chosen
\[
I _{8} = 0  ,\qquad
I_{10} =  {8 \over \delta}
{  (e^{ \delta I_{1} } - 1)^{2}
\over
I_{2} + I_{3} } I_{6},
\]
i.e.
\begin{gather}
{\delta \Delta x + 2 (e^{\delta \tau} -
1)
\left( { h^{-} \over h^{+} + h^{-} } ( \ln u) _{x}  +
    { h^{+} \over h^{+} + h^{-} } ( \ln u )_{\bar{x}}
\right) = 0 } ,\nonumber\\
{\delta (\Delta x) ^{2} + 4 (1 -
e^{-\delta \tau} ) \left( \ln \hat{u} - e^{ \delta
\tau} \ln u \right) = {8 \over \delta} {  (e^{\delta
\tau} - 1)^{2} \over h^{+} + h^{-} }
           \left[ (\ln u) _{x} - ( \ln u)_{\bar{x}}
\right] }  .\label{sh53}
\end{gather}
One invariant solution of this scheme is given
in~\cite{[Bakir]}.

{\bf 4. A linear heat equation without a source
($\boldsymbol{Q = 0}$).}

{ \bf 4.1. Preliminary consideration.}
The linear heat transfer equation
\begin{equation} \label{eq54}
u_{t} = u_{xx}
\end{equation}
admits a six-parameter Lie symmetry group of point
transformations,
corresponding to the infinitesimal operators
\begin{gather}
X_{1} = {{\partial \over \partial t}} ,\qquad X_{2} =
{{\partial \over \partial
x}} , \qquad X_{3} = {2t} {{\partial \over \partial
x}} {-} xu {{\partial \over
\partial u}}  ,\qquad X_{4} = {2t} {{\partial \over \partial t}} + {x}
{{\partial \over \partial x}},\nonumber\\
 X_{5} = {4t^{2}} {{\partial \over \partial t}}
+ {4tx} {{\partial \over
\partial x}} - ( {x^{2}} + {2t} {)u} {{\partial \over
\partial u}} ,
\qquad X_{6} = {u} {{\partial \over \partial u}} \label{op54}
\end{gather}
and an infinite-dimensional symmetry
\[
X^{*} = {a(x,t)}  {{\partial \over \partial u}} ,
\]
where $a(t,x)$ in an arbitrary solution of equation~(\ref{eq54}).
Symmetry $ X^{*} $ represents linearity of the
equation~(\ref{eq54}).

Probably, the simplest approximation of the linear
equation
is the explicit scheme
\begin{equation}  \label{eq55a}
{ \hat{u} - u \over \tau } = { u_{+} -2u + u_{-} \over
 h
^{2} }
\end{equation}
considered on a uniform orthogonal mesh. This equation
is invariant
with respect to the operators $X_1$, $X_2$, $X_4$ and
$X_6$ of the
set~(\ref{op54}). Since the equation is linear
it possesses a~superposition principle that is
reflected in the invariance with respect to the
operator
\[
X^{*}_h  = {a_h (x,t)} {{\partial \over \partial u}},
\]
where  $  {a_h (x,t)} $ is an arbitrary solution of
equation~(\ref{eq55a}).
In~\cite{[Bakir]} it was shown how to construct a
discrete model
which admits the six-dimensional group  (\ref{op54}).

To preserve the Galilean operator $X_3$ and the
projective
operator $X_5$ it is necessary to introduce a moving
mesh.

{\bf 4.2. Heat transfer system of equations
and
superposition principle.}
With the help of the differential invariants of the
operators
(\ref{op54}) in the space $( t, x, u, u_{x}, u_{xx},
dt$, $dx, du)$
\[
J_{1} = { { dx + 2 { { u_{x} \over {u} }
} dt }
\over { dt^{1/2}}} ,
\qquad J_{2} = {{du\over u}} + {{1\over 4}} {{dx\over
dt}}^{2} + \left( -
{{u_{xx}\over u}} + {u_{x}^{2}\over u^{2}} \right) dt
\]
we can represent the heat equation~(\ref{eq54})
as the system
\[
J_{1} = 0, \qquad J_{2} = 0
\]
that is
\begin{equation} \label{sys54}
{dx\over dt}=   -2 {u_{x}\over u },\qquad
{du\over dt} = u_{xx} -
{2}{{u_{x}^{2}\over u}}.
\end{equation}
By construction this system is invariant with respect
to the
six-dimensional group gene\-ra\-ted by the operators~(\ref{op54}).
The system also inherits the superposition principle
of the linear heat transfer equation. The
superposition principle
has the form of  summing two solutions of the system~(\ref{eq54}),
but it also acts on the trajectories on which the
variable~$x$ evolves.
For two arbitrary solutions
$U_{1}(t,x)$ with the trajectories $X_{1}(t)$:
\[
{dX_{1} \over dt } = - { 2U_{1x} \over U_{1} },
\]
and $U_{2}(t,x)$ with the trajectories  $X_{2}(t)$:
\[
{dX_{2} \over dt } = - { 2U_{2x} \over U_{2} },
\]
their linear combination
\begin{equation} \label{sup54}
U= \alpha U_{1} + \beta U_{2},\qquad  \alpha, \beta =
{\rm const},
\end{equation}
is also the solution of (\ref{sys54}). However, this
linear combination has its own trajectories
satisfying
\begin{equation} \label{trac54}
{dX \over dt } = { \alpha U_{1}  \over {\alpha
U_{1}+\beta U_{2} } }
{dX_{1} \over dt } +
                 { \beta U_{2}  \over {\alpha
U_{1}+\beta U_{2} } }
{dX_{2} \over dt }.
\end{equation}
Therefore, the superposition principle can be presented in the
following form:
\begin{gather}
  \left(  \begin{array}{c}
                U(t,x) \vspace{1mm}\\
           { \displaystyle { dX \over dt } }
          \end{array} \right)
=   \left(  \begin{array}{cc}
        \alpha & 0 \vspace{1mm}\\
         0 & { \displaystyle { \alpha U_{1}  \over
{\alpha U_{1}+\beta U_{2} } } }
    \end{array} \right)
\left(  \begin{array}{c}
                U_{1}(t,x) \vspace{1mm}\\
             { \displaystyle { dX_{1} \over dt } }
          \end{array} \right)    \nonumber\\
\phantom{ \left(  \begin{array}{c}
                U(t,x) \vspace{1mm}\\
           { \displaystyle { dX \over dt } }
          \end{array} \right)=}{}+
 \left(  \begin{array}{cc}
           \beta   & 0 \vspace{1mm}\\
         0 & { \displaystyle { \beta U_{2}  \over
{\alpha U_{1}+\beta U_{2} } } }
    \end{array} \right)
\left(  \begin{array}{c}
                U_{2}(t,x) \vspace{1mm}\\
             { \displaystyle { dX_{2} \over dt } }
          \end{array} \right).
\end{gather}

Let us show the superposition principle for system~(\ref{sys54})
by an example. The solution
\[
U_{1} = { 1 \over \sqrt{ t+ t_{1}} } \exp \left( - { (
x - a ) ^{2} \over
4(t+ t_{1} ) } \right)
\]
has the trajectories
\[
x = a + ( x_{0} - a) { t + t_{1} \over t_{1} },
\]
while the solution
\[
U_{2} = { 1 \over \sqrt{ t+ t_{2} } } \exp \left( - {
( x - b ) ^{2} \over
4(t+ t_{2} ) } \right)
\]
exists on the trajectories
\[
x = b + ( x_{0} - b) { t + t_{2} \over t_{2} }.
\]
The linear combination (\ref{sup54}) of these two
solutions is also the solution of
system~(\ref{sys54}). Its trajectories are
\[
{  {
{ dX \over dt} =
{
{
{
{ \alpha  \over \sqrt{ t+ t_{1}} } \exp \left( - { ( x
- a ) ^{2} \over
4(t+ t_{1} ) } \right)
  \left( { x - a   \over t+ t_{1} } \right)
+
{ \beta  \over \sqrt{ t+ t_{2} } } \exp \left( - { ( x
- b ) ^{2} \over
4(t+ t_{2} ) } \right) }
 \left( {  x - b   \over t+ t_{2}  } \right)
}
\over
{
{ \alpha  \over \sqrt{ t+ t_{1}} } \exp \left( - { ( x
- a ) ^{2} \over
4(t+ t_{1} ) } \right) +
{ \beta  \over \sqrt{ t+ t_{2} } } \exp \left( - { ( x
- b ) ^{2} \over
4(t+ t_{2} ) } \right) } }
} } .
\]
Examples of evolution of grid points and corresponding
solutions are shown
on Figs.~5--8 (for computations we used discrete
model (\ref{sh54}) which
will be introduced in point 4.3 of this section).

\begin{figure}[th]
\centerline{\epsfig{figure=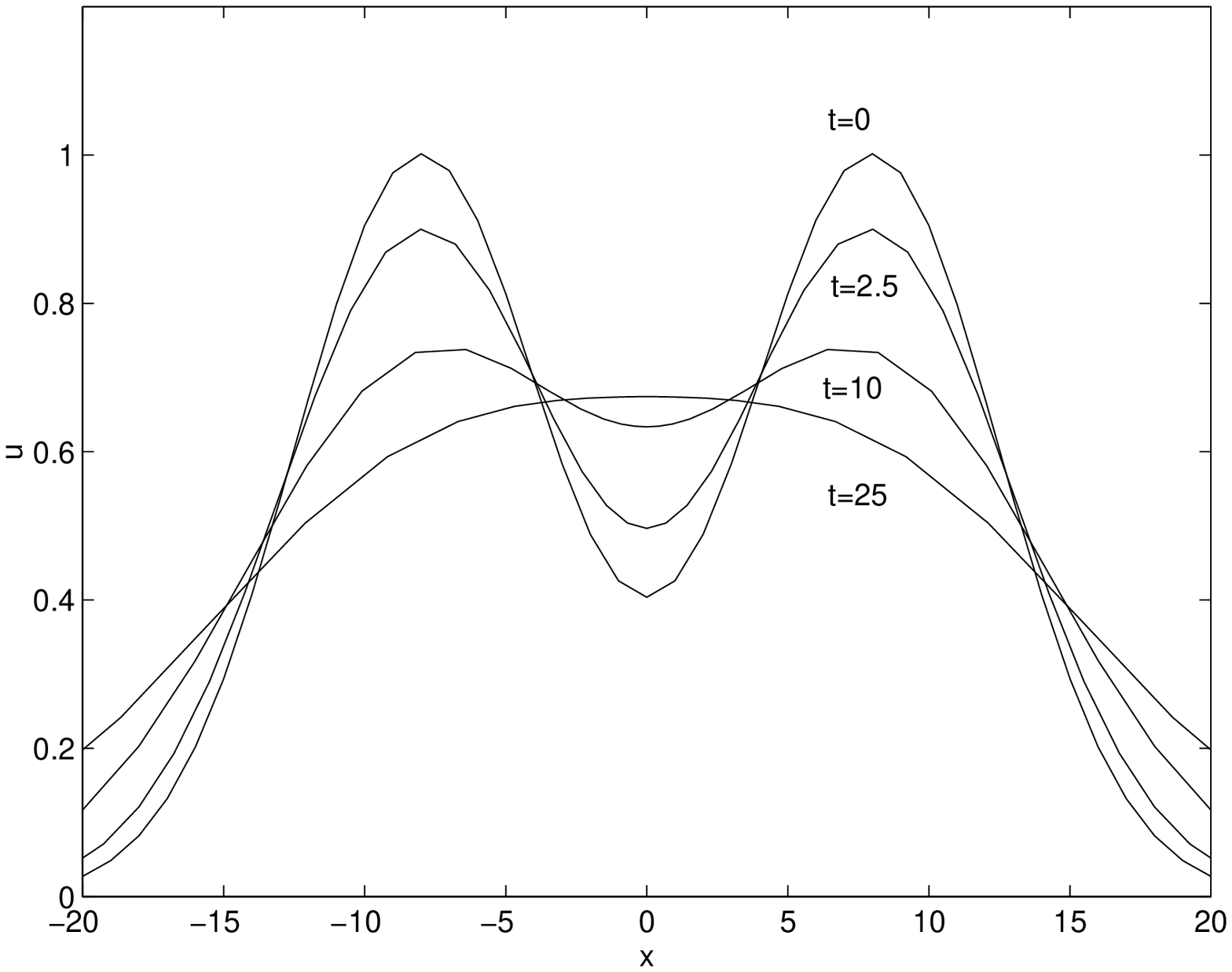,angle=0,width=8cm}}

\vspace{-3mm}

\caption{Evolution of solution (\ref{sup54}):
$\alpha=\beta =1$, $t_1=t_2=10$, $a=-8$, $b=8$. }
\end{figure}
\begin{figure}[th]
\centerline{\epsfig{figure=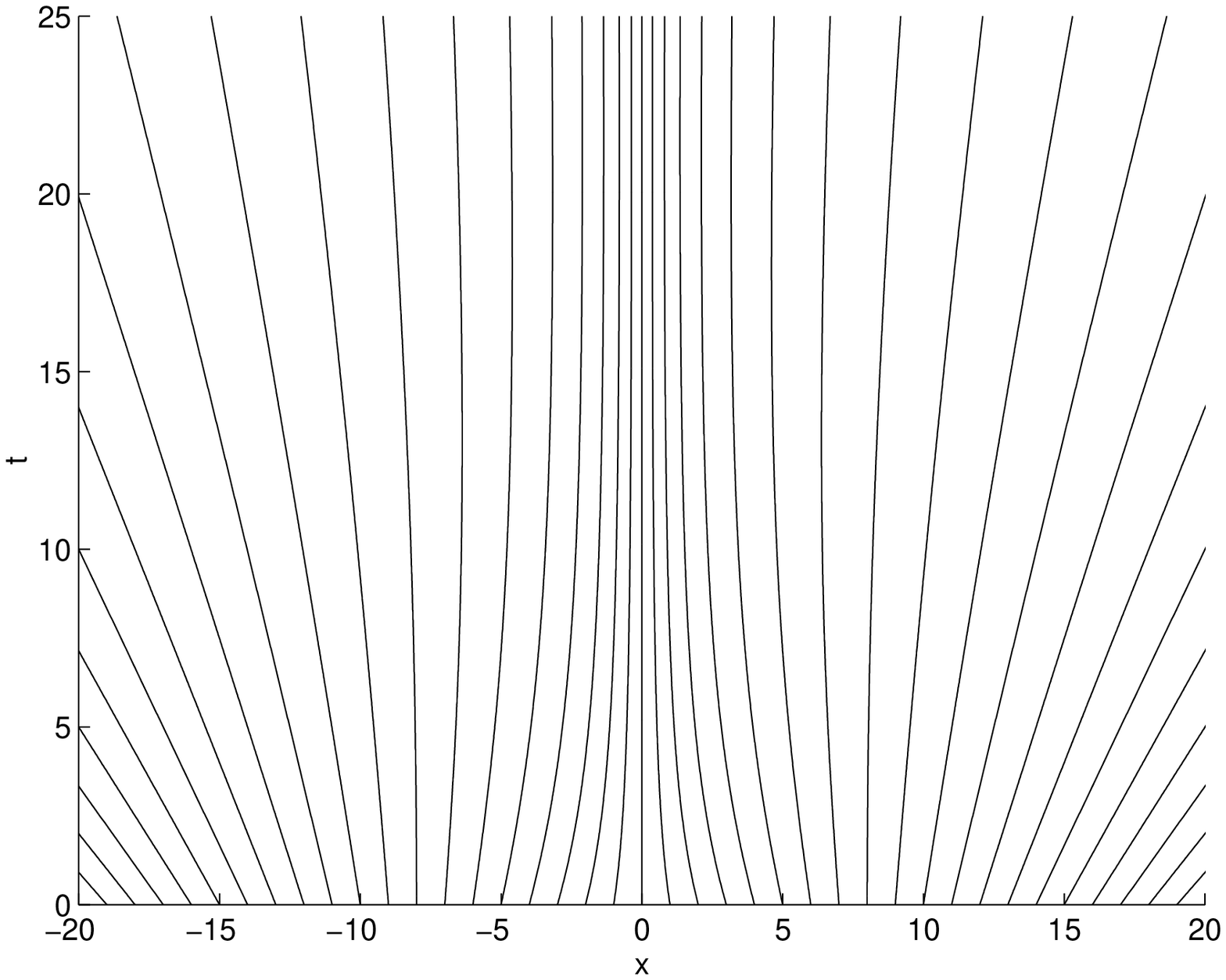,angle=0,width=8cm}}

\vspace{-3mm}

\caption{Mesh trajectories for the solution shown in
Fig.~5.}
\end{figure}
\begin{figure}[th]
\centerline{\epsfig{figure=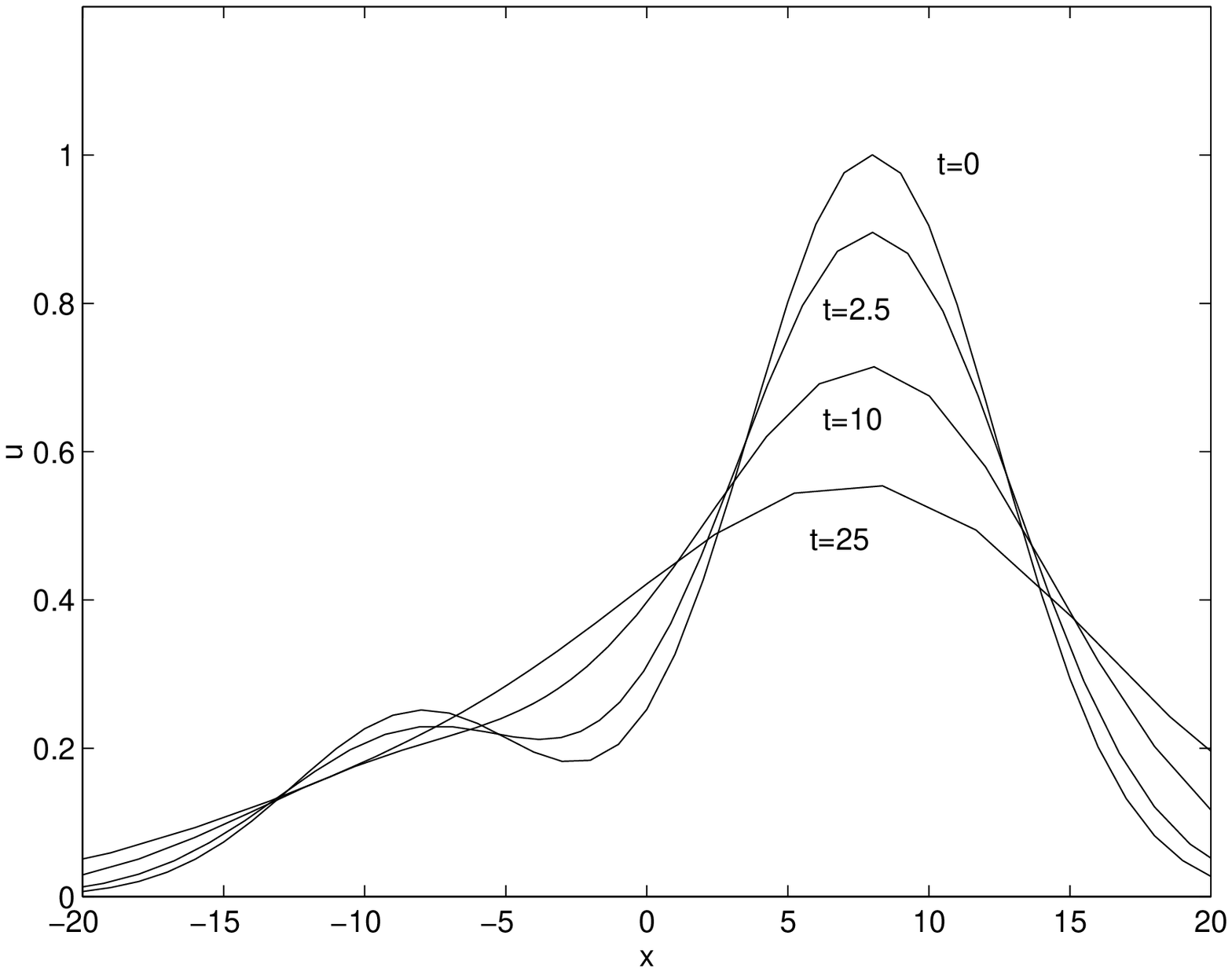,angle=0,width=8cm}}

\vspace{-3mm}

\caption{Evolution of solution (\ref{sup54}):
$\alpha=0.25$, $\beta =1$, $t_1=t_2=10$, $a=-8$,
$b=8$.}
\end{figure}
\begin{figure}[th]
\centerline{\psfig{figure=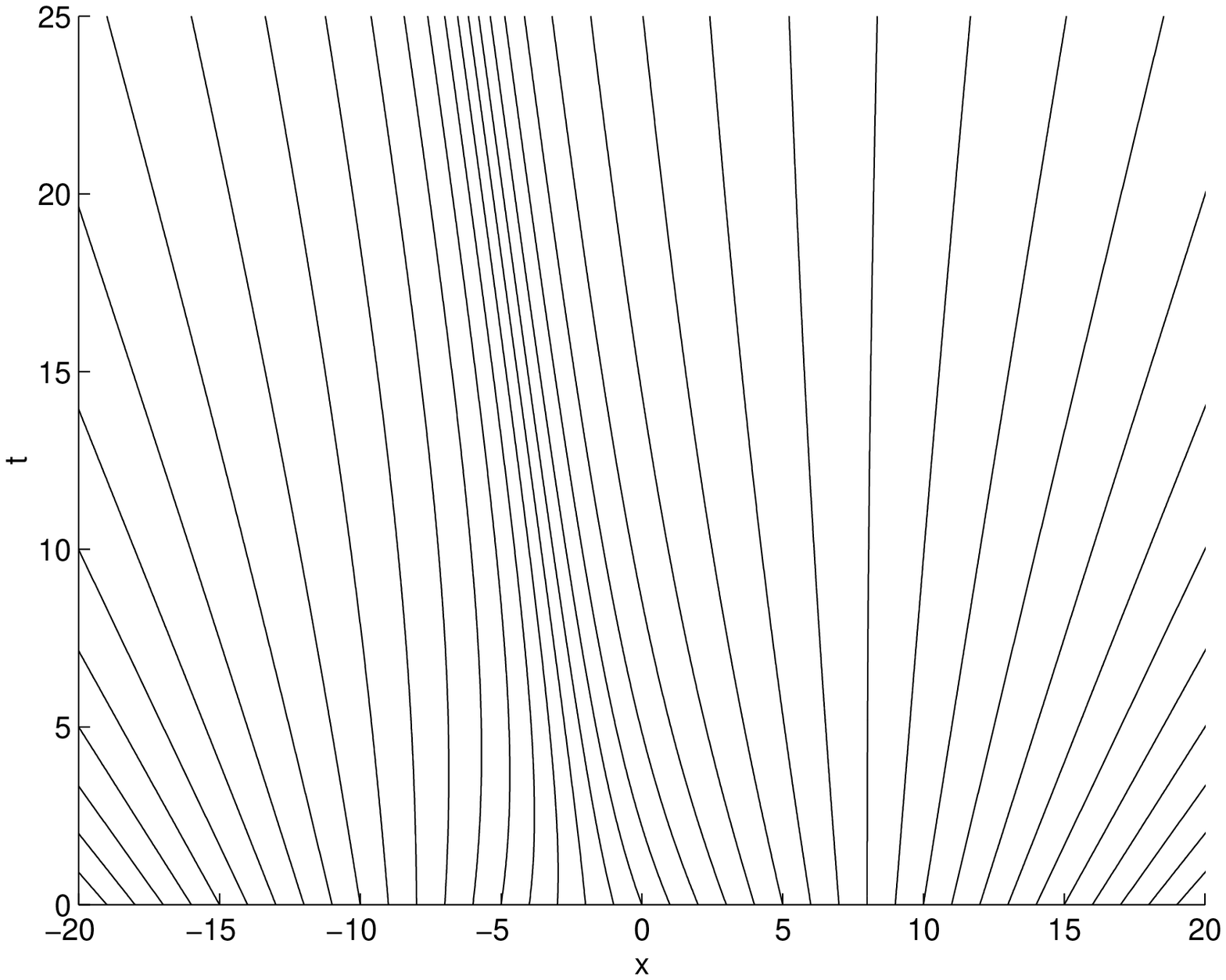,angle=0,width=8cm}}

\vspace{-3mm}

\caption{Mesh trajectories for the solution shown in
Fig.~7.}
\end{figure}

If we consider Lagrangian derivatives of the solution~$U$,
we could get the superposition
principle
\[
\left(  { d U \over dt }
-  U_x   { d X \over dt } \right) =
\alpha
\left(  { d U_{1} \over dt }  -
  U_{1x}   { d X_{1} \over dt }
\right) +
\beta
\left(  { d U_{2} \over dt }  -
  U_{2x}   { d X_{2} \over dt }
\right),
\]
i.e.\ the
superposition principle $U_{t} = \alpha U_{1t} + \beta U_{2t}$
of the linear heat equation expressed in terms of the of total derivatives of~$U$
and~$X$.

{ \bf 4.3. Invariant schemes on moving
meshes.}
For the difference modeling of system~(\ref{sys54}) we
need the whole set
of difference invariants of Lie symmetry group~(\ref{op54}) in the
difference space corresponding to the chosen stencil
$( t, \hat{t}, x,
\hat{x}, h^{+}, h^{-}, \hat{h}^{+}, \hat{h}^{-},u,
\hat{u}, u_{+}$, $u_{-},
\hat{u}_{+}, \hat{u}_{-} )$:
\begin{gather*}
I_{1} = { {h^{+}} \over {h^{-}} } ,
\quad I_{2} = { {\hat h^{+}} \over {\hat h^{-}}} ,
\quad I_{3} = { { {\hat h^{+}} { h^{+}}} \over \tau }
,\quad I_{4} = { \tau^{1/2} \over h^{+} }
{{\hat{u}\over u}} {\exp } \left(
{{1\over 4}} {{(\Delta x)^{2}\over \tau}} \right) ,
\\
I_{5} = {{1\over 4}} {{h^{+2}\over \tau}} -
{{h^{+2}\over h^{+} + h^{-}}}
\left( {{1\over h^{+}}} {\ln }  {{u_{+}\over u}}  +
{{1\over h^{-}}} {\ln }
 {{u_{-}\over u}} \right) ,
\\
I_{6} = {{1\over 4}} {{\hat{h}^{+2} \over \tau}} +
{\hat{h}^{+2} \over
\hat{h}^{+} + \hat{h}^{-}} \left( {1\over \hat{h}^{+}}
{\ln } { \hat u_{+}
\over \hat u} + {1\over \hat{h}^{-}} {\ln } {\hat
u_{-} \over \hat u} \right) ,
\\
I_{7} = {\Delta xh^{+} \over \tau}  + {2h^{+} \over
h^{+} + h^{-}}  \left(
{h^{-} \over h^{+}} {\ln }  {u_{+} \over u}  - {h^{+}
\over h^{-}} {\ln }
 {u_{-} \over u}  \right) ,
\\
I_{8} = {{\Delta x\hat{h}^{+}\over \tau}} +
{{2\hat{h}^{+}\over \hat{h}^{+} +
\hat{h}^{-}}} \left( {{\hat{h}^{-}\over \hat{h}^{+}}}
{\ln } {\hat u_{+} \over
\hat u}  - {{\hat{h}^{+}\over \hat{h}^{-}}} {\ln }
{\hat u_{-} \over \hat u}
\right) .
\end{gather*}

Approximating system (\ref{sys54}) by the invariants,
we obtain a system
of difference evolution equations. As an example, we
present here
an invariant difference model which has explicit
equations for
the solution~$u$ and the trajectory of~$x$:
\begin{gather}
{\Delta x}  =  { {{2\tau\over h^{+}
+ h^{-}}} }\left( -
{{{h^{-}\over h^{+}}} } {\ln }
{{{u_{+}\over u}} }
+ {{{h^{+}\over h^{-}}} } {\ln }
{{{u_{-}\over u}}
} \right),\nonumber\\
{\left(  {{u\over \hat{u}}}  \right)
}^{2} {\exp } \left( - {
{{1\over 2}} {{ (\Delta x)^{2}\over
\tau}} } \right)  =   1
-  {{{4\tau\over h^{+} + h^{-}}} \left(
{{1\over h^{+}}} {\ln }
{{u_{+}\over
u}} + {
{1\over h^{-}}} {\ln } {{u_{-}\over u}} \right) }.\label{sh54}
\end{gather}
We also can write out an implicit model
\begin{gather*}
{\Delta x}  =  { {{2\tau\over
\hat{h}^{+} + \hat{h}^{-}}} \left( -
{{\hat{h}^{-}\over \hat{h}^{+}}} {\ln }
{{\hat{u}_{+}\over \hat{u}}} + {{\hat{h}^{+}\over
\hat{h}^{-}}} {\ln } {{\hat{u}_{-}\over \hat{u}}}
\right) },\\
{\left( {{{\hat{u}\over u}} } \right)
}^{2} {\exp } \left( { {{1\over 2}}
{{\Delta x^{2}\over \tau}} } \right)  =   1 +
{ {{4\tau\over \hat{h}^{+} +
\hat{h}^{-}}} \left( {{1\over \hat{h}^{+}}} {\ln }
{\hat u_{-} \over \hat u} + {{1\over \hat{h}^{-}}}
{\ln } {\hat u_{-} \over \hat u} \right) } .
\end{gather*}
It is also possible to combine an explicit equation
for the mesh
and an implicit approximation of the PDE or vice
versa.
Other ways to approximate the system~(\ref{sys54}) are
also possible.

{\bf 4.4. Optimal system of subalgebras and
reduced systems.}
Among all invariant solutions there is a minimal set
of such solutions,
called the optimal system of invariant solutions. 
From this set of invariant solutions
any invariant solution can be obtained by an
appropriate group
transformation. 
The difference model (\ref{sh54}) is a system of two
evolution
equations. To find its invariant solutions we need to
provide a time mesh which is invariant with respect to
the considered operator. An invariant time mesh giving
flat time layers
can be represented by the equation
\begin{equation} \label{meshtt}
\tau_i = g ( t_i )  ,\qquad i = 0,1,2,\ldots.
\end{equation}
We request this equation to be invariant with respect
to
the considered symmetry. Since for the operators~(\ref{op54}) the
coefficients $ \xi^t$ do not depend on~$x$ and~$u$ we
can propose
an invariant time mesh for any symmetry.
In the case   $ \xi^t = 0 $ the function~$g$ can be
taken arbitrary. For example, we can choose the uniform mesh
$ t_j = j \tau $, $ \tau = {\rm const}$.
Thus, different invariant solutions may have different time meshes.

The adjoint action of the Lie group transforms an
invariant solution
into another one~\cite{[16],[15]}. In our  case
it also transforms the time mesh equation~(\ref{meshtt}). Thus the
adjoint action gives us a new invariant solution with
a~corresponding invariant mesh.

On  the example of the difference model (\ref{sh54})
we will construct the optimal system of solutions
which are invariant with respect to one-parameter
groups. The optimal system of one-dimensional subalgebras of
the algebra of symmetries for the linear heat equation consists of
algebras corresponding to the operators (see~\cite{[15]}):
\begin{gather*}
Y_{1}=X_{2}={\partial \over \partial x}, \qquad
Y_{2}=X_{6}=u{\partial \over \partial u}, \qquad
Y_{3}=X_{1}+cX_{6}={ \partial \over \partial t}
+ cu {\partial \over \partial u},
\\
Y_{4}=X_{1} - X_{3}= {\partial \over \partial t}
- 2t { \partial \over \partial x}
+ xu {\partial \over \partial u}, \qquad
Y_{5}=X_{4}+2cX_{6}= 2t {\partial \over \partial t}
+ x {\partial \over \partial x}
+ 2cu {\partial \over \partial u},
\\
Y_{6}=X_{1}+X_{5}+cX_{6}=(4t^{2} +1) {\partial \over
\partial t}
+ 4tx {\partial \over \partial x}
+ (c-x^{2}-2t)u{\partial \over \partial u}.
\end{gather*}

Let us find invariant solutions corresponding to these
one-dimensional subalgebras.

1) The subalgebra corresponding to the operator
$Y_{1}$ has only constant solutions $ u = C $, $ C =
{\rm const}$
considered on an orthogonal mesh $ \Delta x = 0 $.

 2) The subalgebra corresponding to the
operator $Y_{2}$ does not
have invariant solutions (the necessary condition of
existence of invariant
solutions does not hold~\cite{[16]}).

 3) The operator $Y_{3}$ has the following
invariants:  $u\exp(-ct)$,
$\tau$ and $\Delta x$. The time step~$\tau$ is
invariant, so we can consider
a uniform time mesh. We will seek a solution of the
difference
model in the form
\[
u=\exp(ct)f(x).
\]
Substituting this invariant form of the solution into
system~(\ref{sh54}), we get:
\begin{gather}
\Delta x = {{ 2\tau \over h^{+} + h^{-}
} } \left( - { { h^{+} \over h^{-} } }
\ln \left( { { f(x+h^{+}) \over f(x) } }
\right) + { { h^{+} \over h^{-} } } \ln
\left( { { f(x-h^{-}) \over f(x) } }
\right) \right),\nonumber\\
\left( {{ f(x) \over f(x+ \Delta x) } }
\right)^{2}       \exp \left( -2c \tau -
{{1\over 2} } {{ \Delta
x^{2} \over \tau} } \right) \nonumber\\
\qquad{}
=  1 -{ { 4\tau \over h^{+} + h^{-} }
} \left( { { 1\over h^{+}} } \ln \left( {
 { f(x+h^{+}) \over f(x)} } \right) + {
{ 1 \over h^{-}} } \ln \left( {
 { f(x-h^{-}) \over f(x) } } \right)
\right) .
 \label{opt3}
\end{gather}

System (\ref{opt3}) will become  a system of  two ordinary difference
equations if we project it into the invariants space.
To project the system we have to impose
\begin{equation} \label{opt3a}
\Delta x    =
-h^{-}, \ 0, \  \mbox{or} \  h^{+}.
\end{equation}
A solution of system (\ref{opt3}) with one of
conditions~(\ref{opt3a})
provides the solution of system~(\ref{sh54}) which is
invariant with respect to the operator~$Y_{3}$.

 4) The operator $Y_{4}$ has invariants: $u
\exp \left(-xt -{2\over3}t^{3} \right)$, $x+t^{2}$,
$\tau$ and ${\Delta x \over 2\tau} - t$.
Let us search a solution of the difference model
(\ref{sh54}) in the form
\[
u = \exp \left(tx +{2\over 3}t^{3} \right) f(x+t^{2}).
\]
By means of variables
\begin{gather*}
y=x+t^{2},\qquad y-h^{-}_{y}=x-h^{-}+t^{2},
\\
y+h^{+}_{y}=x+h^{+} + t^{2},\qquad y + \Delta y = x +
\Delta x + ( t + \tau )^{2}
\end{gather*}
we get the following system for
the invariant solution of system (\ref{sh54}):
\begin{gather*}
\Delta y - \tau^{2}  =  { { 2\tau
\over h^{+}_{y} + h^{-}_{y} } } \left( - {
 { h^{-}_{y} \over  h^{+}_{y} } } \ln
\left( {{ f(y+h^{+}_{y}) \over f(y) } }
\right) +  { { h^{+}_{y} \over h^{-}_{y}
} } \ln \left( { { f(y- h^{-}_{y})
\over f(y) } } \right) \right),\\
\left( { { f(y) \over f(y+ \Delta y)} }
\right)^{2} \exp \left( - { { 1\over 2
\tau} } \Delta y^{2}  - \tau (2y + \Delta y) +
{ { 1\over 6} }\tau^{3} \right)  \\
\qquad{}=  1 - { {4\tau \over h^{+}_{y} +
h^{-}_{y} } } \left( {  { 1\over
h^{+}_{y}}  }  \ln \left( {
{f(y+h^{+}_{y})\over f(y)} } \right) + {{1\over h^{-}_{y}} } \ln \left( {  {
f(y-h^{-}_{y} \over f(y) } } \right) \right),
\end{gather*}
where $\Delta y$ can have one of the following values
\begin{equation} \label{opt4a}
\Delta y   =
-h^{-}_{y}, \  0 \  \mbox{or} \  h^{+}_{y}.
\end{equation}
A solution of the above system with one of conditions~(\ref{opt4a})
let us find the invariant solution for the operator~$Y_{4}$.

 5) Expressions
${x \over \sqrt t}$ , $t^{-c}u$, ${\tau \over t}$ and
${\Delta x \over x}$
are invariants of the operator $Y_{5}$.
Let us search a~solution of the difference model in
the form
\[
u= t^{c}f\left({x \over \sqrt{t} } \right).
\]
In variables
\[
y = { x \over \sqrt{\displaystyle
t}},\qquad
y - h_{y}^{-} ={ x - h^{-} \over
\sqrt{ t}},
\qquad
y + h_{y}^{+} = {  x +
h^{+} \over \sqrt{t}},\qquad y +
\Delta y = {x + \Delta x \over
\sqrt{ t+ \tau}}
\]
we get the following system of equations:
\begin{gather*}
\sqrt{ {1+a}} ( y + \Delta y ) - y  =
{  {2a \over h^{+}_{y} + h^{-}_{y} }
}
\left( - {{h^{-}_{y} \over h^{+}_{y} }
} \ln
\left( { { f(y+h^{+}_{y}) \over f(y) } }
\right)
+ {{ h^{+}_{y} \over h^{-}_{y} } } \ln
\left( { { f(y-h^{-}_{y}) \over f(y) }
}\right) \right), \\
({1+a})^{-2c}
\left( { {f(y) \over f(y+ \Delta y)} }
\right)^{2}  \exp
\left( - { {1\over 2} } \left( (y+ \Delta
y)
\sqrt { \displaystyle {1+a \over a} }
- y { { 1\over \sqrt{ \displaystyle a}}
}\right)^{2} \right)  \\
\qquad{}=  1 - { { 4a \over h^{+}_{y} +
h^{-}_{y} } } \left( { { 1 \over
h^{+}_{y}} } \ln \left( { {
f(y+h^{+}_{y}) \over f(y) } } \right) + {
 { 1\over h^{-}_{y} } } \ln \left( {
 { f( y-h^{-}_{y} ) \over f(y) } }
\right) \right).
\end{gather*}
Here $\Delta y$ can have one of the values determined
by conditions~(\ref{opt4a}) and $a$ is a constant from the condition
$a={\tau \over t}$
which determines an invariant time spacing.
This condition can be found if we look for
a time spacing $\tau = g(t)$ which is invariant with
respect to the
considered operator~$Y_{5}$.

6) For the operator $Y_{6}$ we have the
following invariants:
\begin{gather*}
{x \over \sqrt{4t^{2}+1}}, \qquad (4t^{2}+1)^{1/4}u
\exp \left( {tx^{2} \over 4t^{2}+1} + {c\over2}
\arctan(2t) \right),
\\
{ 4t^{2}+ 1 \over \tau} + 4 t, \qquad {\Delta x \over
x} {4t^{2} + 1 \over \tau} - 4t.
\end{gather*}
We look for a solution of the difference model in the
form
\[
u = (4t^{2}+1)^{-1/4} \exp \left( - {tx^{2} \over
4t^{2}+1} - {c\over 2} \arctan(2t) \right) f \left( {
x \over \sqrt {4t^{2}+1}} \right).
\]
Involving new variables
\begin{gather*}
y={\displaystyle  x\over \sqrt{4t^{2}+1}},\qquad y - h^{-}_{y} = {  x -
h^{-} \over \sqrt{ 4t^{2}+1}},\\
y + h^{+}_{y} = {x + h^{+} \over
\sqrt{4t^{2} +1}},\qquad y + \Delta y =
{x + \Delta x \over \sqrt{4(t+\tau)^{2}+1}},
\end{gather*}
the system of equations~(\ref{sh54}) can be presented
in the form:
\begin{gather*}
\sqrt{{b^{2} + 1} } ( y + \Delta y) - b
y   =
 {{ 1\over h^{+}_{y} + h^{-}_{y} } }
\left(  \!- {{ h^{-}_{y} \over h ^{+}_{y}
} } \ln \left( {{ f( y+h^{+}_{y} )
\over f(y) } } \right) \!+ { { h^{+}_{y}
\over h_{y}^{-} } } \ln \left( {{
f(y-h_{y}^{-}) \over f(y)} } \right) \right),\\
\sqrt{ { b^{2}+1} } \left(
{{ f(y)\over f(y+ \Delta y)} }
\right)^{2}\\
\qquad {}\times
\exp \left( c \arctan \left( {  {1 \over
b} }\right) - b (y^{2}+(y+ \Delta y)^{2}) + 2 \sqrt{
{ b^{2}+1} } y(y+\Delta y) \right) \\
\qquad {}=  b - {{ 2\over h^{+}_{y} +
h^{-}_{y} } } \left( { { 1 \over
h^{+}_{y} }  }  \ln \left( { {
f(y+h^{+}_{y}) \over f(y)} } \right) + {{ 1\over h^{-}_{y} } } \ln \left( {  {
f(y-h^{-}_{y} ) \over f(y) } } \right) \right).
\end{gather*}
where $\Delta y$ has one of the values of~(\ref{opt4a}),
$b$ is the  constant from necessary condition of
invariant grid existence
\[
2b= 4t + {4t^{2}+1 \over \tau }.
\]

Therefore, the obtained reduced systems of equations
determine the optimal system of invariant solutions
for the difference model of the liner heat transfer
equation. It means that each invariant solution can be
found by transformation
of a solution from the optimal system with the help of 
the corresponding element of the group.
As we mentioned before the invariant time mesh for the
new solution is obtained from the time mesh of the
solution from the
optimal system with the help of the same group
transformation.
For example, the transformation
corresponding to the operator $X_{1}$
with the value of the parameter $-t_{0}$ gives shift
in time
$\hat{t} = t - t_{0}$. Since
\[
\mbox{Ad}\, [\exp(-t_{0}X_{1})]  Y_{5} = X_{*} = Y_{5}
+ 2t_{0}X_{1},
\]
the action of this transformation transfers the
invariant solution
with respect to the operator $Y_{5}$ into the
solution which is invariant with respect to the
operator $X_{*}$.
By this transformation the
spacing ${\tau \over t}=a$ is transformed into the
spacing
${\tau \over t + t_{0} } = a$.

 {\bf Example of an exact solution.}
Among all group invariant solutions for the difference
model~(\ref{sh54}) there is one interesting solution which can be
integrated exactly~\cite{[Bakir]}.
This is the solution invariant with respect to the
operator
\begin{equation} \label{oper54}
2  t_0  X_{2} + X_{3}, \qquad  t_0 = {\rm const},
\end{equation}
namely the solution
\begin{equation} \label{eqr57}
u(x,t) = C \left( { t_0 \over t + t_0 } \right) ^{1/2}
\exp \left( - {x^{2} \over 4(t + t_0 )} \right),
\end{equation}
considered on the mesh
\begin{equation} \label{eqr557}
x_i = x^{0}_i \left( {t + t_0 \over t_0 } \right) ,
\end{equation}
where $ x^{0}_i$ are space mesh point at $ t = t_0$.
In the case $ t_0 = 0 $ we get the well known
fundamental solution
of the linear heat equation. Note that it has a
``singular'' mesh.

Let us show how this solution can be obtained form
the optimal system of the invariant solutions. From
\[
\mbox{Ad}\, [ \exp( \varepsilon X_{5}) ] Y_{1} =
X_{2} +  2 \varepsilon    X_{3}
\]
we see that the solution (\ref{eqr57}) can be obtained
from the
solution invariant with respect to operator $Y_{1}$
by the transformation
$ \mbox{Ad}\, [ \exp( \varepsilon X_{5}) ] $ with $
\varepsilon = { 1 \over 4 t_0} $.
If we take the original solution on the orthogonal
mesh which is
uniform in space and has the following special time
spacing on the interval
$[ 0 , t_0]$:
\[
u_i^j = C , \qquad x_i = i h   , \qquad  i = 0, \pm 1,
\pm 2, \ldots,   \qquad
t_j = {   j \tau  t_0 \over  t_0 +  j \tau } i , \qquad
 j = 0,1,2,\ldots
\]
then the proposed transformation provides us the
solution~(\ref{eqr57}) on the uniform space mesh~(\ref{eqr557})
and uniform time mesh~$ t_j = j \tau$.

Thus, we see that difference model (\ref{sh54})
inherits both the group admitted by the original differential
equation and the ability to be integrated on a subgroup.

{\bf 4.5. The way to stop a moving mesh.}
The obtained difference models have adap\-ti\-ve
nonorthogonal grids.
We can find a way to stop the moving mesh, i.e., an
exchange of variables
which orthogonalizes the mesh.
The differentiation operator of Lagrange type~${d\over
dt}$
can be presented in the following form
\[
{d \over dt} = D_{t}  - 2 {u_{x} \over u} D_{x} ,
\]
where
\[
D_{t} = {\partial \over {\partial t}} + u_{t}
{\partial \over {\partial u}}
+ {\cdots } ,
\qquad
D_{x} = {\partial \over {\partial x}} + u_{x}
{\partial \over {\partial u}}
+ {\cdots } .
\]
The operator ${d\over dt}$ in contrast to the
operators $D_{t}$ and $ D_{x}$
does not commute with the operators of total
differentiation with respect to~$t$ and~$x$:
\[
\left[ {d\over dt}  , D_{t}  \right] =
2 \left( { u_{xt} \over u} - { u_{x} u_{t} \over u^{2}
} \right) D_{x} ,
\qquad
\left[ {d\over dt} , D_{x} \right] =
2 \left( { u_{xx} \over u} - { u_{x} ^{2} \over u^{2}
} \right) D_{x} .
\]
It is necessary to find an operator of total
differentiation with respect to
a new space variable $s$ such that
\begin{equation} \label{opdif}
\left[  {d \over dt }  ,  D_{s} \right] = 0.
\end{equation}
The last commuting property is possible if we involve
a new dependent
variable $ \rho > 0 $ (density)~\cite{[17]}.
The operator $D_{s}= {1\over \rho } D_{x}$  satisfies
(\ref{opdif}) if
$\rho$ holds the equation
\[
\rho _{t} - 2\rho  \left( {u_{xx}\over u} -
{u^{2}_{x}\over u^{2}} \right) - 2 {u_{x}\over u} \rho
_{x} = 0 .
\]
The new space variable $s$ is introduced with the help
of equations
\[
s_{t} =  { {2\rho } {u_{x} \over u}
}, \qquad
s_{x} = \rho.
\]
For convenience we can put initial data $ \rho ( 0 ,
x) \equiv 1$.
Then, $ s = x$ for $ t = 0$.

In the variables $(t,s)$ the heat transfer equation
will get a form of the system
\begin{equation} \label{ts1}
{u_t} =  {\rho ^{2}}
\left( u_{ss} - 2{u^{2}_{s} \over u} \right) + {\rho
\rho _{s} u_{s}} ,\qquad
{ \rho_t} = {2\rho
^{3}} \left( {u_{ss} \over u} - {u^{2}_{s} \over
u^{2}} \right) + {2\rho ^{2}\rho _{s}} {u_{s} \over u}
\end{equation}
which can be rewritten in the form of conservation
laws
\begin{equation} \label{ts2}
\left( {1\over \rho } \right) _t  =
\left( -  {2\rho } {u_{s} \over u} \right) _s ,\qquad
\left( {u\over \rho } \right) _t   =
( - \rho u_{s} ) _s .
\end{equation}
The space coordinate $x$
is defined by the system of equations
\begin{equation} \label{ts3}
x_{t} = - 2 \rho  {u_{s}\over u} ,\qquad
x_{s} = {1\over \rho }  .
\end{equation}
System (\ref{ts1}) in the space of independent
variables $(t,s)$ and
extended set of dependent variables $(u,\rho,x)$
admits a group of point transformations determined by
the following infinitesimal operators:
\begin{gather}
X_{1} = {\partial \over \partial t} , \qquad X_{2} =
{\partial \over \partial x} , \qquad X_{3} = 2t
{\partial \over \partial x} - xu {{\partial \over
\partial u}} ,
\nonumber\\
X_{4} = 2t {{\partial \over \partial t}} + x
{{\partial \over \partial x}} + s {{\partial \over
\partial s}} , \quad X_{5} = 4t^{2} {{\partial \over
\partial t}} + {4tx} {{\partial \over \partial x}} - (
x^{2} + 2t )u {\partial \over \partial u} - {4t\rho }
{{\partial \over \partial \rho }} ,
\nonumber\\
X_{6} = u {{\partial \over \partial u}} , \qquad X^{*}
= {f(s)} {\partial \over \partial s} + {\rho
f^{\prime}(s)} {{\partial \over \partial \rho }} ,\label{ts4}
\end{gather}
where $f(s)$ is an arbitrary function of $s$.

In the independent variables $(t,s)$ operators
$X_{1}$--$X_{6}$ are operators of
Lie algebra factorized by the operator $X^{*}$.
Condition~(\ref{cht}) of grid orthogonality and
condition~(\ref{ch}) of
space grid uniformity hold and it gives an opportunity
to construct
a difference model which is invariant
with respect to operators $X_{1}$--$X_{6}$ on the
orthogonal grid.

Let us write system (\ref{ts1}), (\ref{ts3}) in the
form of differential invariants. In the space of
variables $(t, x, s, u, \rho , dt, dx, ds, du,
d\rho , u_{s}, \rho _{s},x_{s},   u_{ss} )$
there  are five invariants:
\begin{gather*}
J_{1} = x_{s} \rho  , \qquad
J_{2} = {\rho \over ds} \left(dx + {2\rho } {u_{s} \over u}
dt \right) , \qquad
J_{3} = {(ds)^{2} \over \rho ^{2} dt} ,
\\
J_{4} = {(ds)^{2} \over \rho ^{3}} \left( {d\rho \over
dt} - {\rho _{s}ds \over dt} - 2\rho ^{3} \left(
{u_{ss} \over u} - {u^{2}_{s} \over u^{2}} \right)
-2\rho ^{2} \rho _{s} {u_{s} \over u} \right) ,
\\
J_{5} = \left( {ds \over \rho } \right)^{2} \left( -
{2\over u} {du \over dt} - {1\over 2} \left( {dx\over
dt} \right)^{2} + 2\rho ^{2} \left( {u_{ss} \over u} -
{u^{2}_{s} \over u^{2}} \right)
+ 2\rho \rho _{s} { u_{s} \over u} \right) .
\end{gather*}
With the help of these invariants we rewrite system
(\ref{ts1}), (\ref{ts3}) as
\begin{gather}
{{u_t} } = { \rho^{2}
\left( u_{ss} - {u_{s}^{2} \over u} \right) + \rho
\rho_{s} u_{s} },\nonumber\\
{ \rho_t } = {2
\rho^{3} \left( {u_{ss} \over u} - {u^{2}_{s} \over
u^{2}} \right) + 2 \rho^{2} \rho_{s} {u_{s} \over u}
}, \qquad
{ x_t} = -2 \rho {{u_{s} \over u} }  \label{ts5}
\end{gather}
with the constraint equation $ x_s = {1\over \rho} $.

Now we can find a system of equations which
approximates~(\ref{ts5})
and is invariant with respect to the set of operators~(\ref{ts4}).
We can use the six-point invariant stencil (Fig.~4).
We have the following invariants for the set of
operators~(\ref{ts4}) where
the operator~$X^{*}$ is changed on its difference
analog
\[
X^{*}_h  = f(s) {\partial \over \partial s} + \rho  {
\dsd D } ( f(s) )
{\partial \over \partial \rho}
\]
in the corresponding to the chosen stencil
space
$( t, \hat{t}, s,
h^{+}_{s},h^{-}_{s}, x, \hat{x}, h^{+}_{x}, h^{-}_{x},
\hat{h}^{+}_{x}, \hat{h}^{-}_{x}, u, \hat{u}$, $u_{+},
 u_{-}, \hat{u}_{+},
\hat{u}_{-}, \rho , \hat{\rho }, \rho _{+},
\rho _{-}, \hat{\rho }_{+}, \hat{\rho }_{-})$:
\begin{gather*}
I_{1} = {h^{+}_{x} \over h^{-}_{x}} , \quad
I_{2} = {\hat{h}^{+}_{x} \over \hat{h}^{-}_{x}}  ,
\quad
I_{3} = { h^{+}_{x} \hat{h}^{+}_{x} \over \tau }
,\quad
I_{4} = { \tau^{1/2} \over h^{+}_{x}}  {\hat{u} \over
u} {\exp } \left( {1\over 4} {\Delta x^{2} \over \tau}
\right) ,
\\
I_{5} = {1\over 4} {h^{+2}_{x} \over \tau } -
{h^{+2}_{x} \over {h^{+} + h^{-}}}  \left( {1\over
h^{+}_{x}} {\ln}  {u_{+} \over u} + {1\over h^{-}_{x}}
{\ln}  {u_{-} \over u} \right) ,
\\
I_{6} = {1\over 4} {\hat{h}^{+2}_{x} \over \tau } +
{\hat{h}^{+2}_{x} \over \hat{h}^{+}_{x} +
\hat{h}^{-}_{x}} \left( {1\over \hat{h}^{+}_{x}} {\ln}
 {\hat{u}_{+} \over \hat{u}} + {1\over
\hat{h}^{-}_{x}} {\ln } {\hat{u}_{-} \over \hat{u}}
\right) ,
\\
I_{7} = {\Delta x h^{+}_{x} \over \tau} + {2h^{+}_{x}
\over h^{+}_{x} + h^{-}_{x}} \left( {h^{-}_{x} \over
h^{+}_{x} } {\ln } {u_{+} \over u} - {h^{+}_{x} \over
h^{-}_{x}}  {\ln } {u_{-} \over u} \right) ,
\\
I_{8} = {\Delta x \hat{h}^{+}_{x} \over \tau} +
{2\hat{h}^{+}_{x} \over \hat{h}^{+}_{x} +
\hat{h}^{-}_{x} }   \left( {\hat{h}^{-}_{x} \over
\hat{h}^{+}_{x}}  {\ln } {\hat{u}_{+} \over \hat{u}} -
 {\hat{h}^{+}_{x} \over \hat{h}^{-}_{x}} {\ln }
{\hat{u}_{-} \over \hat{u}} \right) ,
\\
I_{9} = {\hat{\rho }_{-} \over \rho _{-}} , \quad
I_{10} = {\hat{\rho} \over \rho } , \quad
I_{11} = {\hat{\rho }_{+} \over \rho _{+}} ,  \quad
I_{12} = {h^{+}_{s} \over \rho  h^{+}_{x}} , \quad
I_{13} = {h^{-}_{s} \over \rho _{-} h^{-}_{x}} .
\end{gather*}

With the help of these invariants we can write
difference model in the form
of the following system of evolution difference
equations
(we present here only one invariant difference model
which corresponds to the system (\ref{sh54}) in
variables $(t,x)$):
\begin{gather*}
{\Delta x}  =  { {2\tau} {
{ { - { h^{-}_{s} \over h^{+}_{s} } {
\rho \over \rho_{-} } \ln { u_{+} \over u } +
{h^{+}_{s} \over h^{-}_{s} } { \rho_{-} \over \rho }
\ln { u_{-} \over u }} }  \over  { {
h^{+}_{s} \over \rho } + { h^{-}_{s} \over \rho _{-}}
} } },\\
\left( { {u\over \hat{u}} } \right) ^{2}
{\exp } \left( {- {1\over 2} {\Delta
x^{2}\over \tau} } \right)  =  1 -
{ {4\tau} { { {\rho \over
h^{+}_{s}} \ln  { u_{+} \over u } +  { \rho_{-} \over
h^{-}_{s} }  \ln {  u_{-} \over u } } \over
{ { h^{+}_{s}  \over \rho } + { h^{-}_{s}
\over \rho _{-} } } } },\\
{\hat{\rho}   =  \rho   { h^{+}_{x}  \over
\hat{h}^{+}_{x} } } .
\end{gather*}
In the case of the uniform grid ($h^{+}_{s}=
h^{-}_{s} = h_{s}$)
this model could be simplified as follows:
\begin{gather*}
\Delta x  =  { {2\tau}  {
{ - \rho^{2} \ln { u_{+} \over u }  +
\rho _{-} ^{2} \ln { u_{-} \over u } } \over h_{s} (
\rho + \rho_{-}) } },\\
\left( { {u\over \hat{u}} } \right)^{2}
{\exp } \left( - { {1\over 2} {\Delta
x^{2}\over \tau} } \right)  =   1 -
{ {4\tau \rho \rho _{-} \over h^{2}_{s} (
\rho  + \rho _{-} )}  \left( {\rho } {\ln } {u_{+}
\over u} +  {\rho _{-}} {\ln } {u_{-}\over u} \right)
},\\
{\hat{\rho}   =  \rho
{ h^{+}_{x}  \over \hat{h}^{+}_{x} }}  .
\end{gather*}

The system (\ref{ts1}) has only two dependent
variables $u$ and $ \rho$
and it can be approximated without involvement of the
space variable~$x$.
However, Galilean symmetry  $X_3$ and projective
symmetry~$X_5$ are nonlocal in the coordinate system~$(t,s)$
and we need to consider the dependent variable~$x$
in order to have these symmetries. Constructing an
invariant with respect to the set of operators~(\ref{ts4})
difference model, we inevitably involve~$x$ into the
difference equations.

It is important to notice that in all cases moving
meshes could be stopped by using the Lagrange type coordinate
system (for an introduction of Lagrange type coordinate systems see, for example,~\cite{[17]}).

{\bf 5.} If $Q=\delta u$, $\delta = \pm 1$, the
equation
\begin{equation} \label{e55}
u_{t}=u_{xx} + \delta u
\end{equation}
can be transformed into equation (\ref{eq54}) by the
change of variables
\begin{equation} \label{ch55}
\bar{u} = u e^{ - \delta t }.
\end{equation}
Reversing this transformation, one can get an
invariant model for equation~(\ref{e55})
from an invariant model for the heat transfer equation
without a source.

{\bf 6.} $Q=\delta = {\rm const}$. The equation has the
form:
\begin{equation} \label{e56}
u_{t}=u_{xx} + \delta.
\end{equation}
The case $Q=0$, which was considered in point 4 of
this section,
is a partial case for $Q=\delta = {\rm const}$, but the
constant source
can be excluded by the evident transformation
\begin{equation} \label{ch56}
\bar{u} = u - \delta t.
\end{equation}
It means that we can get a difference model for~(\ref{e56}) from the model~for (\ref{eq54}).

\section{Concluding remarks}

In the paper [10] the entire set of invariant schemes for ordinary
difference equations of the second order was developed. There were
shown that for some equations and symmetries it is necessary to involve
nontrivial lattices, which are not uniform in a space of independent variable
and should depend on solution. Such schemes are self adapted for any solution
and they are as much exactly integrable, as its continuous counterpart.

In the papers [2, 3, 8, 9, 11] several examples of invariant schemes for PDEs
with two independent variables (KdV equation, nonlinear Schr\"odinger
equation etc.) were constructed. Again, there were shown that for
some equations and symmetries it is necessary to involve nontrivial
two-dimensional meshes, which are not uniform and rectangular in a
space of two independent variables and should depend on solution. Such schemes
are self-adapted for any solution, the meshes are evolutionary in time and these
schemes have as many exact invariant solutions as their continuous counterparts.
Moreover, for invariant variational cases invariant schemes have difference
conservation laws as well as for continuous case.

Thus, the applications of symmetry to difference equations led to the
evolution of idea of possible meshes: from simple regular stationary
meshes to self-adapted moving meshes, explicitly depending on solutions.
Notice, that this idea evolution well corresponds to the big changes
in numerical analysis, where the idea of self-adaptivity of schemes
and meshes is in a broad fashion now.

In this paper we developed the entire set of invariant difference
schemes for the heat transfer equation
\begin{equation} \label{equ}
u_t = \left( K(u) u_{x} \right) _{x} + Q(u),
\end{equation}
for arbitrary coefficients $K(u)$,
$Q(u)$ and for all special cases of the coefficients which extend the
symmetry group admitted by equation (\ref{equ}).

The main conclusion is that we have presented an algorithmic
way to construct the invariant difference schemes (i.e.\ a difference
equation and a mesh it is defined on) for all cases of underlying heat equation.

Other conclusion is that symmetry preservation in difference
schemes led to essential different discrete models. For different cases
of coefficients $K(u)$ and~$Q(u)$ taken in accordance with the group
classification of equation~(\ref{equ}) we have obtained
different discrete models: for some cases $K(u)$
and $Q(u)$ we had to construct discrete models not for
equation (\ref{equ}), but for the equivalent Lagrangian system:
\begin{equation} \label{syst}
{dx\over dt} = \varphi(u,u_{x}),
\qquad
 {du\over dt}= \psi(u,u_{x},u_{xx}).
\end{equation}
The consideration of the equation (\ref{equ}) in the form of
system~(\ref{syst}) at the very beginning
would provide us with the classification by functions~$\varphi$ and~$\psi$.
In that case there is one-to-one correspondence between the
coordinate systems for continuous and discrete cases
of the group classification.
In particular, if the symmetry of the equation~(\ref{equ})
does not require application of Lagrangian type moving mesh
(evolving in time mesh), then we have $\varphi \equiv 0$. In that case
the classification of the system~(\ref{syst}) and corresponding
difference equations on orthogonal mesh can be carried out by means of~$\psi(u,u_{x},u_{xx})$.

\subsection*{Acknowledgments}

This work was sponsored in part by Russian Fund for
Base Research and The Norwegian Research Council under contracts
no.111038/410,  through the SYNODE project, and no.135420/431, through
the BeMatA program.

\label{dorodnitsyn-lastpage}
\end{document}